\documentclass[generic,preprint]{imsart}
\RequirePackage{natbib}

\startlocaldefs

\usepackage{amssymb,amsmath,amsfonts,amsthm,amscd,amstext,mathtools}
\usepackage[utf8]{inputenc}
\usepackage{xspace,enumerate}
\usepackage{graphicx}
\usepackage[pdftex]{hyperref}
\usepackage{color}
\usepackage[charter]{mathdesign}
\usepackage{subcaption}

\usepackage[a4paper,scale={0.72,0.74},marginratio={1:1},footskip=7mm,headsep=10mm]{geometry}

\numberwithin{equation}{section}

\newtheorem{theorem}{Theorem}[section]
\newtheorem{claim}[theorem]{Claim}
\newtheorem{lemma}[theorem]{Lemma}
\newtheorem{proposition}[theorem]{Proposition}
\newtheorem{corollary}[theorem]{Corollary}
\newtheorem{remark}[theorem]{Remark}
\newtheorem{definition}[theorem]{Definition}

\long\def\xcom#1{}

\newcommand{\cA}{{\ensuremath{\mathcal A}} }

\newcommand{\cC}{{\ensuremath{\mathcal C}} }
\newcommand{\cD}{{\ensuremath{\mathcal D}} }
\newcommand{\cE}{{\ensuremath{\mathcal E}} }

\newcommand{\cL}{{\ensuremath{\mathcal L}} }

\newcommand{\cV}{{\ensuremath{\mathcal V}} }

\newcommand{\gep}{\varepsilon}       % \ge already exists...

\newcommand{\tx}[1]{\texttt{#1}}

\renewcommand{\tilde}{\widetilde}          % wider `tilde'
\DeclareMathSymbol{\leqslant}{\mathalpha}{AMSa}{"36} % nicer `smaller or equal'
\DeclareMathSymbol{\geqslant}{\mathalpha}{AMSa}{"3E} % nicer `larger or equal'
\DeclareMathSymbol{\eset}{\mathalpha}{AMSb}{"3F}     % nicer `emptyset'
\newcommand{\R}{\mathbb{R}}

\newcommand{\Z}{\mathbb{Z}}
\newcommand{\N}{\mathbb{N}}

\newcommand{\ind}{{\sf 1}}

\renewcommand{\epsilon}{\varepsilon}
\DeclareMathOperator\argmax{arg\, max}

\newenvironment{myitemize}{%
\begin{list}{$\bullet$}%
        {%
        \setlength{\itemsep}{0.4em}%
        \setlength{\topsep}{0.5em}%
        \setlength\leftmargin{2.45em}%
        \setlength\labelwidth{2.05em}%
        \setlength{\labelsep}{0.4em}%
        }%
        }%
{\end{list}}

\renewenvironment{itemize}{
\begin{myitemize}}%
{\end{myitemize}}

 \newcommand{\be}[1]{\begin{equation}\label{#1}}
 \newcommand{\ee}{\end{equation}}

 \newcommand{\bl}[1]{\begin{lemma}\label{#1}}
 \newcommand{\el}{\end{lemma}}

 \newcommand{\br}[1]{\begin{remark}\label{#1}}
 \newcommand{\er}{\end{remark}}

 \newcommand{\bt}[1]{\begin{theorem}\label{#1}}
 \newcommand{\et}{\end{theorem}}

 \newcommand{\bd}[1]{\begin{definition}\label{#1}}
 \newcommand{\ed}{\end{definition}}

 \newcommand{\bcl}[1]{\begin{claim}\label{#1}}
 \newcommand{\ecl}{\end{claim}}

 \newcommand{\bp}[1]{\begin{proposition}\label{#1}}
 \newcommand{\ep}{\end{proposition}}

 \newcommand{\bc}[1]{\begin{corollary}\label{#1}}
 \newcommand{\ec}{\end{corollary}}

 \newcommand{\bpr}{\begin{proof}}
 \newcommand{\epr}{\end{proof}}

 \newcommand{\bi}{\begin{itemize}}
 \newcommand{\ei}{\end{itemize}}

\newcommand{\ens}[1]{\left\{#1\right\}}

\newcommand{\valabs}[1]{\left|#1 \right|}

\newcommand{\unsur}[1]{\frac{1}{#1}}

\date{\today}

\endlocaldefs
\begin{document}

\begin{frontmatter}
\title{Interacting partially directed self-avoiding walk: a probabilistic perspective}
\runtitle{Interacting partially directed self-avoiding walk}

\begin{abstract} \footnote{\today}
We review some recent results obtained in the framework of the 
2-dimensional Interacting Self-Avoiding Walk (ISAW). After a brief presentation of the rigorous results that have been obtained so far 
for ISAW we focus on the \emph{Interacting Partially Directed Self-Avoiding Walk} (IPDSAW), a model introduced in 
\cite{ZL68} to decrease the mathematical complexity of ISAW. 

In the first part of the paper, we discuss how a new probabilistic approach based on a random walk representation (see \cite{NGP13}) allowed for a sharp determination of the  asymptotics of the free energy close to criticality (see \cite{CNGP13}).  
Some scaling limits of IPDSAW were conjectured in the physics literature (see e.g. \cite{POBG93}).  We discuss here the fact 
that all limits are now proven rigorously, i.e., for the extended  regime in \cite{CP15}, for the  collapsed regime in \cite{CNGP13} and  
at criticality in \cite{CarPet17a}. 

The second part of the paper starts with the description of four open questions related to physically relevant extensions of IPDSAW. 
Among such extensions is the Interacting Prudent Self-Avoiding Walk (IPSAW) whose configurations are those of the 2-dimensional prudent walk. We discuss the main results obtained in  \cite{PT17} about IPSAW and in particular the fact that its collapse transition is proven to exist rigorously.

\end{abstract}

\author{\fnms{Philippe}
  \snm{Carmona}\corref{}\ead[label=e2]{philippe.carmona@univ-nantes.fr}
\ead[label=u1,url]{http://www.math.sciences.univ-nantes.fr/~carmona/}}

%\and
\author{\fnms{Gia Bao} \snm{Nguyen}\ead[label=e3]{ nguyengb@kth.se}}

%\and
\author{\fnms{Nicolas} \snm{Pétrélis}\ead[label=e4]{nicolas.petrelis@univ-nantes.fr}}
\affiliation{Université de Nantes}

%\and
\author{\fnms{Niccolò}
  \snm{Torri}\corref{}\ead[label=e1]{niccolo.torri@upmc.fr}
\ead[label=u1,url]{http://www.lpma-paris.fr/pageperso/torri/}}

\affiliation{Université de Nantes}

\address{Laboratoire de Math\'ematiques Jean Leray UMR 6629\\
Universit\'e de Nantes, 2 Rue de la Houssini\`ere\\
BP 92208, F-44322 Nantes Cedex 03, France\\ \printead{e2}\\\printead{e4}
}

\address{Department of Mathematics\\
KTH Royal Institute of Technology\\
Lindstedtsv\"agen 25, SE-100 44 Stockholm, Sweden\\ \printead{e3}
}

%\address{KTH - Department of Mathematics\\
%                               Lindstedtsv\?agen 25\\
%                               10044 Stockholm - Sweden\\ \printead{e3}}
%

\address{CNRS \& Université Pierre et Marie Curie,\\
Laboratoire de Probabilités et Modèles Aléatoires,\\
Campus de Jussieu, 4 place Jussieu 75252 Paris Cedex 5\\ \printead{e1}\\%\printead{e1}
}

\runauthor{P. Carmona et al.}



\thanks{{\it Acknowledgements.}  P. Carmona, N. P\'etr\'elis and N. Torri thank the Centre Henri Lebesgue ANR-11-LABX-0020-01 for creating an attractive mathematical environment.}

\end{frontmatter}

\vspace{1cm}
\noindent
* Invited paper to appear in a special volume of the Journal of Physics A, on the occasion of the 75th birthday of Stu Whittington.

\tableofcontents

%%%%%%%%%%%%%%%%%%%%%%%%%%% SECTION 1 %%%%%%%%%%%%%%%%%%%%%%%%%%%%%%%%%%%%%%%%%%%%

\section{Introduction}
\label{S1}

The collapse transition is a well known example of phase transition. It   
takes place for instance when an homopolymer is dipped in a poor solvent. 
As the solvent temperature decreases, it reaches a threshold (the $\theta$-point) below which the geometry of a typical polymer configuration changes drastically so that it looks pretty much like a compact ball.

A good mathematical model to investigate this phenomenon is the Interacting Self-Avoiding Walk (see \cite{TF94} or \cite{S86}).
In size $L\in \N$, the configurations of ISAW are given by the $L$-step self-avoiding walk trajectories on $\Z^d$.  A Gibbsian weight is assigned to each such configuration as  $\beta\in [0,\infty)$ (the interaction intensity) times the number of {\it self-touchings}, i.e., pairs of sites of the walk adjacent on the lattice though not consecutive along the walk.  Among lattice polymer models, the ISAW plays a central role because it fulfills the excluded volume effect, a feature that real world polymers indeed satisfy. However, few mathematical results are available so far, mostly because the mathematical understanding of self-avoiding walks remains fairly incomplete. At the moment, the existence of the free energy is established for small interaction parameter $\beta$ (first in  \cite{Uelt02} for random walk with infinite range step distribution and more recently in \cite{HH17} for a larger class of a priori laws on the walk including the simple random walk) but remains open elsewhere.  In dimension $d\geq 5$ and for small 
$\beta$, the mean square displacement of ISAW is proven to be of order $L$  (see \cite{Uelt02})   by using lace expansion techniques.
There is so far, for $d\geq 2$, no rigorous proof of the existence of a phase transition for ISAW.

 \begin{figure}
      %  \centering
        \begin{subfigure}[b]{0.52\textwidth}
                \includegraphics[width=\textwidth]{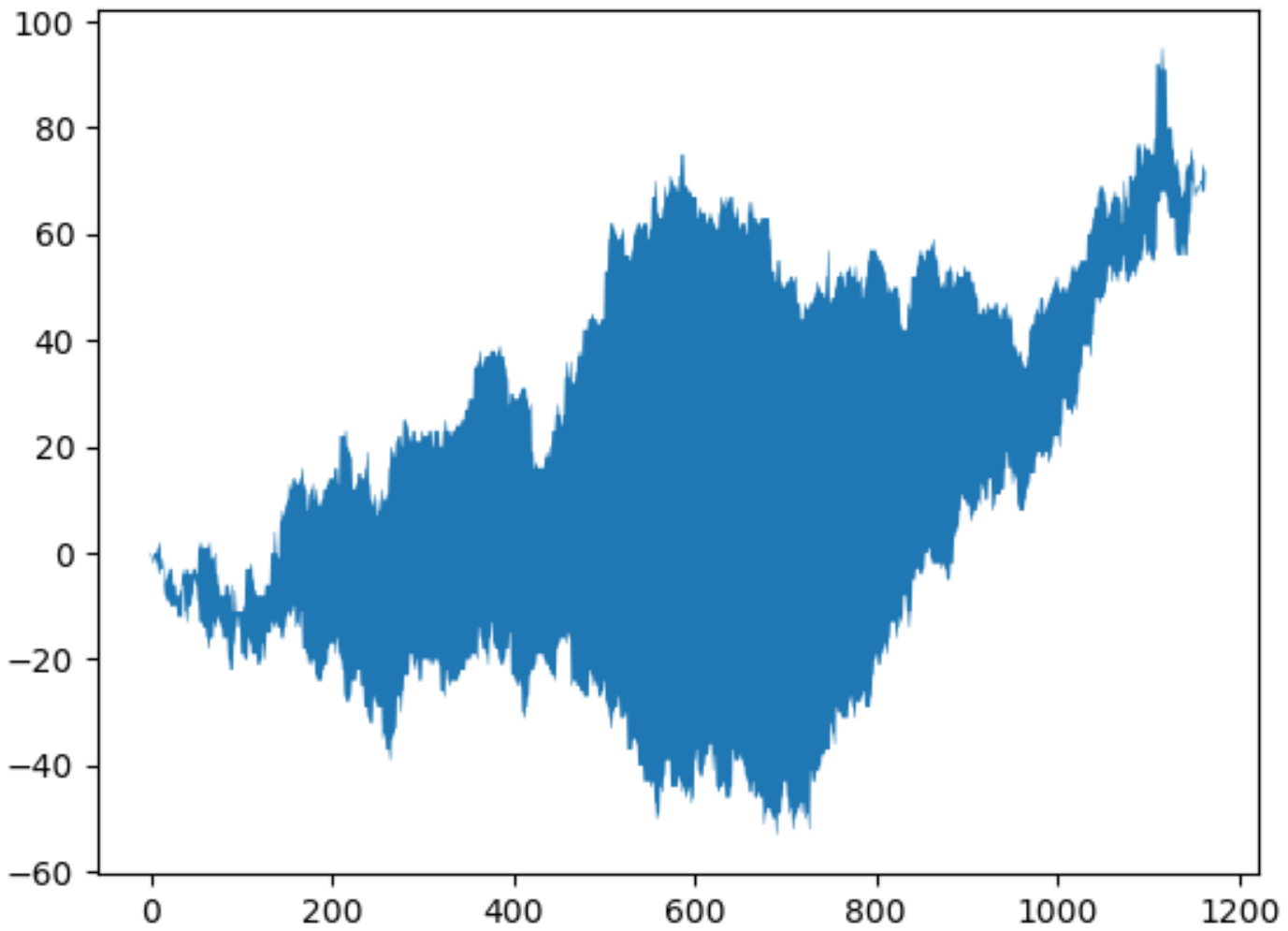}
                \vspace{-3.5cm}
        \end{subfigure}%
	\begin{subfigure}[b]{0.52\textwidth}
\includegraphics[width=\textwidth]{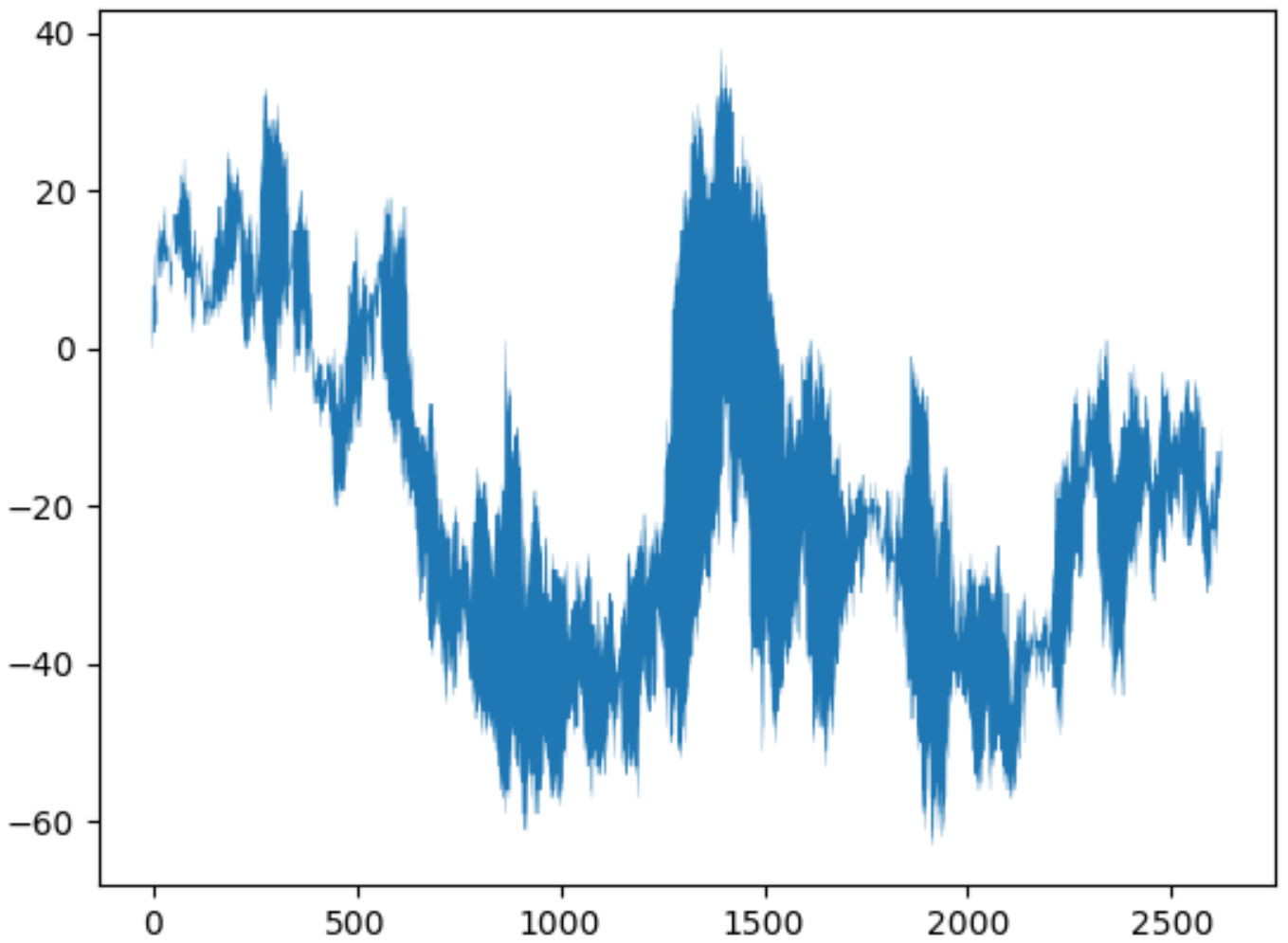}
  \vspace{-3.5cm}
\end{subfigure}
	\caption{Two exact simulations of  IPDSAW at critical temperature
          $\beta=\beta_c$ and with length $L=60000$}
	\label{fig:simulations}
\end{figure}

The mathematical complexity of ISAW has motivated the 
introduction of alternative models for self-interacting random walk. The challenge consists in 
designing models that, on one hand, are  sophisticated enough  to capture the most important physical features 
of the collapse phenomenon and, on the other hand, are tractable enough to allow for a deep mathematical investigation. 
In the physics literature, a lot of attention has been dedicated to exactly solvable models. For instance in \cite{DupSal87}
a two-dimensional polymer model is investigated on the honeycomb lattice. A random environment is introduced by deleting some faces of the lattice in a percolation-type fashion. The edges of the missing faces are prohibited so that, by annealing on the environnement, the resulting model displays attractions between edges. The collapse transition of the model occurs when the deleted faces start to percolate and thanks to this analogy the critical exponents could be computed.  Recent works support the idea that such exactly solvable models
share common features with ISAW itself at criticality. In this spirit, numerical evidences are displayed in \cite{G13} to illustrate the correspondence between the two-dimensional ISAW at criticality and $\text{SLE}_6$ and both theoretical and numerical results are displayed in 
\cite{VJS15} to try and determine the $\theta$-point of ISAW (we actually  refer to the introduction in  \cite{VJS15} for a concise and very clear state of the art on such exactly solvable models).  
Let us now focus on the mathematics literature where two other variants of ISAW received most of the attention.
 
The first of these variants is the Interacting Weakly-Self-Avoiding Walk (IWSAW), introduced in \cite{vdHK01}. In size $L\in \N$,  the set of allowed configurations for IWSAW is much larger than that of ISAW since it contains every $L$-step simple random walk trajectory on $Z^d$. 
However, the Hamiltonian of an IWSAW trajectory contains an additional term that penalizes the auto-contact, i.e., 
decreases the Gibbs weight by $-\gamma$ for every self-intersection of the trajectory. 
The phase diagram of IWSAW is conjectured to be divided into three phases, 
i.e., localized, collapsed and extended. In \cite{vdHK01}, a critical curve $\beta=2d \gamma$ is proven to separate 
the localized phase ($\beta>2d\gamma$) inside which a typical trajectory remains in a box of finite size from the rest of the quadrant. 
Another critical curve $\gamma\mapsto \beta_c(\gamma)$ is conjectured to exist inside $\{(\gamma,\beta)\in [0,\infty)^2\colon \beta<2d\gamma\}$ that separates  a collapsed phase where the end to end distance of a typical trajectory should be $L^{1/d}$ from an extended phase where this distance should be of the same order as that of the self-avoiding walk. In the limit $\gamma\to \infty$, it is expected that 
$\beta_c(\gamma)$ converges to the $\theta$-point of ISAW. Recently, a continuous time version of IWSAW was investigated in \cite{BSW17}. In dimension $4$, an area of the quadrant is isolated (corresponding to small $\gamma$ and $\beta$) and proven to be part of the extended phase.

The second variant is the Interacting Partially-Directed Self-Avoiding walk (IPDSAW) and was first introduced in \cite{ZL68}. This is a 2-dimensional model where 
the set of allowed configurations is narrowed (compared to that of ISAW) but the Hamiltonian remains unchanged. Until recently (see Section \ref{secIPRSAW}) the IPDSAW was the only 2-dimensional polymer model for which the  
collapse transition was rigorously established.  It was first studied with transfer matrix methods (see \cite{BOVY90}) and then
with combinatorial tools in \cite{BGW92} to compute the critical point $\beta_c$ that partitions the phase diagram into 
an extended phase $\cE:=[0,\beta_c)$ and  a collapsed phase $\cC:=[\beta_c, \infty)$.

A  new probabilistic approach of IPDSAW has been introduced in \cite{NGP13} which turned out to strongly simplify its investigation.
In the present paper we review the results obtained using this new framework concerning the analytic properties of free energy in \cite{CNGP13} and \cite{PT17} and also concerning the path properties of IPDSAW in \cite{CP15} and \cite{CarPet17b}.  For every result that we state here, we provide a sketch of its rigorous proof.

%
%\cite{Uelt02}, 
%
%This is the reason why, in the mathematical literature, collapse transition models were rather built by 
%either relaxing the self-avoiding feature of the paths (see for instance or 
%\cite{vdHKK02}) or by considering partially directed paths. 
%This is the case for  the \emph{interacting partially directed self-avoiding walk} (referred to as IPDSAW) 
%that was introduced in \cite{ZL68} and subsequently studied 
%in e.g.  \cite{BGW92} or   \cite{GP13}, \cite{CGP13}
%and \cite{CP15}).
%
%

%%%%%%%%%%%%%%%%%%%%%%%%%%%%%%%%%%%%%%%%%%%%%%%%%%%%%%%%%%%%%%%%%%%%%%%%
%\begin{figure}
%\begin{center}
%\includegraphics[scale = 0.3]{emulsions2.eps}
%\end{center}
%\caption{An undirected copolymer in an emulsion.}
%\label{fig-PolInRa}
%\end{figure}
%%%%%%%%%%%%%%%%%%%%%%%%%%%%%%%%%%%%%%%%%%%%%%%%%%%%%%%%%%%%%%%%%%%%%%%%%

\subsection{The model IPDSAW}
\label{S1.1}

\subsubsection{Mathematical description of the model}
The IPDSAW can be defined in two equivalent manners. In the original definition (see \cite{ZL68}), the polymer configurations 
are modeled by the trajectories of a two dimensional self-avoiding walk, taking unitary steps up, down and to the right,
whereas in the alternative definition, the configurations are modeled by families of vertical stretches. 
In Sections \ref{free} and \ref{Geocha}, we will use the second definition since it fits with the probabilistic approach that we wish to display. 
However, in Section \ref{secIPRSAW} we will come back to the original definition to present the IPSAW, a non-directed extension of
IPDSAW that has been investigated in \cite{PT17}.

\begin{figure}[ht]\center
	\includegraphics[width=.32\textwidth]{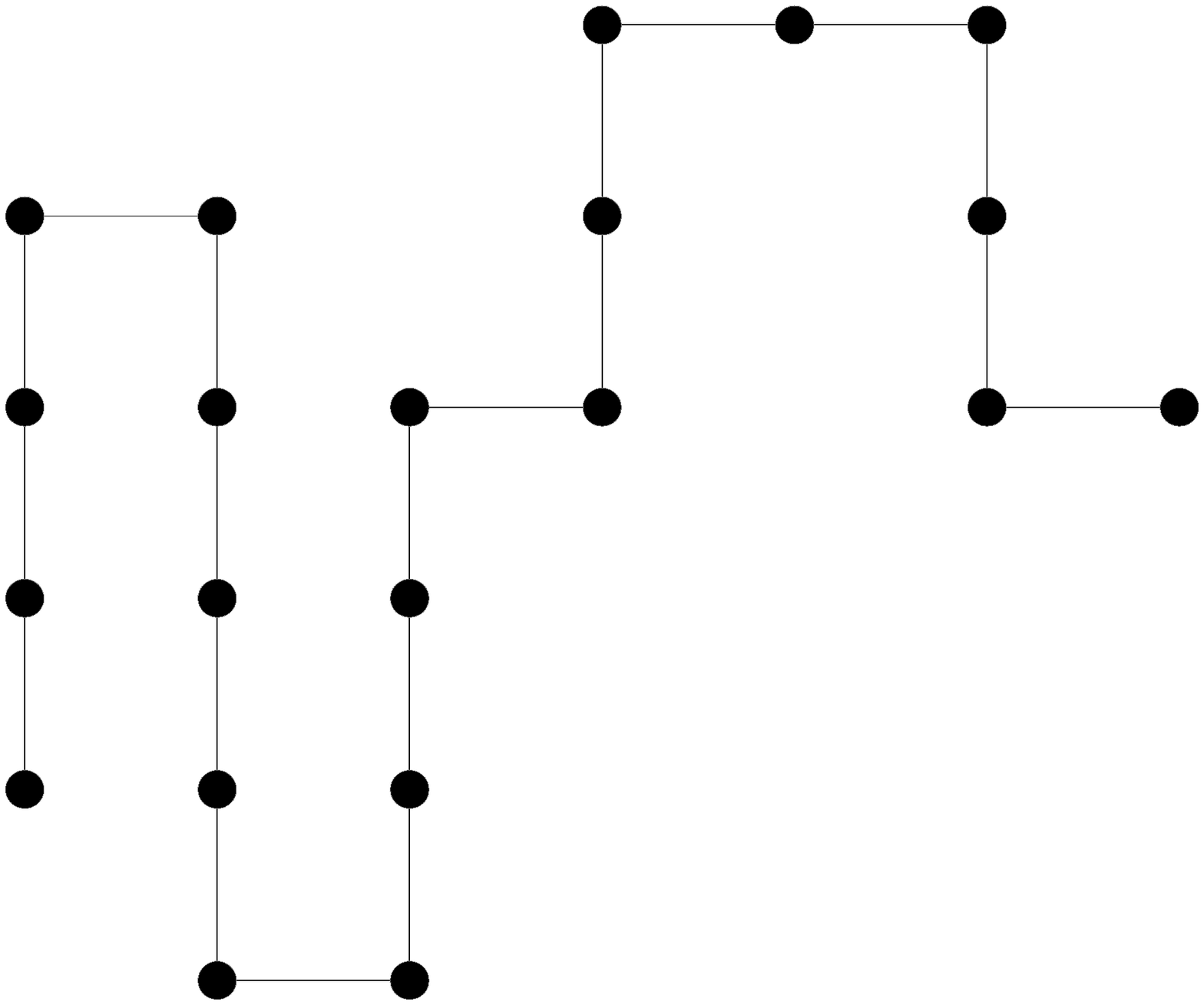}
	\vspace{-1cm}
	\caption{Example of a trajectory  $l\in \cL_{N,L}$ with $N=6$ vertical stretches, a total  length $L=20$ and 
	an Hamiltonian $H_{L}(l)=6 $.}
	\label{fig:stretches}
\end{figure}

In size $L$, the allowed configurations of the polymer can be represented as families of oriented vertical stretches, i.e, 
 $\Omega_L:=\bigcup_{N=1}^L\mathcal{L}_{N,L}$, with
\begin{equation}\label{defLL}
\textstyle\mathcal{L}_{N,L}=\bigl\{l\in\mathbb{Z}^N:\sum_{n=1}^N|l_n|+N=L\bigr\}.
\end{equation}
With such configurations, the modulus of a given stretch
corresponds to the number of monomers constituting this stretch and two consecutive vertical stretches 
are separated by one horizontal monomer  (see Fig. \ref{fig:stretches}).
The repulsion exerted by the solvent on the monomers is taken into account by 
assigning to every configuration $l\in \Omega_L$  an energetic reward $\beta\in [0,\infty)$ every times it performs a \emph{self-touching} 
that is every time it places two non-consecutive monomers at distance $1$ from each other. 
By summing those microscopic interactions, we obtain for $N\in \{1,\dots,L\}$ the Hamiltonian of a given configuration 
$l\in \mathcal{L}_{N,L}$ as  
\begin{equation}
\textstyle H_{L}(l_1,\ldots,l_N)=\sum_{n=1}^{N-1}(l_n\;\tilde{\wedge}\;l_{n+1}),
\end{equation}
where
\begin{equation}
x\, \tilde\wedge\, y=\begin{dcases*}
	|x|\wedge|y| & if $xy<0$,\\
  0 & otherwise.
  \end{dcases*}
\end{equation}
The preceding Hamiltonian is an exponential Gibbs weight that allows us to define the polymer measure on $\Omega_L$ as 
\be{defpolme}
P_{L,\beta}(l)=\frac{e^{\beta H_{L}(l)}}{Z_{L,\beta}},\quad l\in \Omega_L,
\ee
where $Z_{L,\beta}$ is  the partition function of the model. Finally,  the free energy 
\be{deff}
f(\beta):=\lim_{L\to \infty} \frac{1}{L} \log Z_{L,\beta}
\ee
 provides us with the exponential growth rate of the partition function.
 \medskip

 \subsection{Challenges}\label{challenges}

We can distinguish between two main types of questions  that one tries to address when investigating lPDSAW:
\begin{enumerate}
\item \emph{Determine the asymptotic development of the free energy close to criticality.}
\rm We will see below that the free energy of  IPDSAW is trivial in its collapsed phase, i.e., 
$f(\beta)=\beta$ when $\beta\geq \beta_c$. Therefore, one wants to exhibit $\gamma,\alpha >0$ such that 
\be{asydev}
\tilde f(\beta_c-\gep)= \alpha \gep^\gamma (1+o(1))\quad \text{as}\quad  \gep\to 0+,
\ee
with $\tilde f(\beta):=f(\beta)-\beta$ the \emph{excess free energy} of the system. 
One also expects that $\alpha$ may be expressed as the free energy of a counterpart 
model built with Brownian trajectories.
\smallskip

\item \emph{Display the scaling limit of IPDSAW in each regime.}
Compute the growth speed of the horizontal and vertical extensions of a typical 
IPDSAW trajectory in the extended phase ($\beta<\beta_c$), inside the collapsed phase ($\beta>\beta_c$) 
and at criticality $(\beta=\beta_c)$. With those typical growth speeds in hand, determine the 
limiting shape of an appropriately rescaled typical trajectory of IPDSAW. 
\smallskip

%\item \emph{Enrich the model.} Build an advanced version of IPDSAW by enlarging the set of allowed trajectories. 
%\textcolor{blue}{ In this spirit, building for instance a non-directed model of a self-interacting random walk and prove that it undergoes a collapse transition.
%Introducing a non-directed model, for which we can prove a collapse transition,
%interpolating between IPDSAW and ISAW, for which the collapse transition is only conjectured (see for instance \cite{S86})
% represents a step forward in the comprehension of the self-interacting model built with general \emph{self-avoiding} paths.
%% As far as I know, this is the first time that the existence of a critical point is proved for such non-directed model: 
%}
\end{enumerate}

Section \ref{free} below is dedicated to issues of type 1. 
With  Section \ref{Geocha} we settle entirely the issues of type 2. In Section \ref{open} we list some open problems related to ISAW and with Section \ref{secIPRSAW} we give a first answer to one of them. %\textcolor{blue}{discussing the issues of type 3}.

\vspace{-1.1cm}

\section{Asymptotics of the free energy close to criticality}\label{free}
 \subsection{A probabilistic representation of the partition function}\label{stillfree}

In  \cite{NGP13}, 
a Random Walk representation of IPDSAW has been introduced 
(see Section \ref{Geocha} below
for more details). With this new technique, a probabilistic expression of the partition function has been derived, i.e., 
\begin{align}\label{tgh}
\tilde Z_{L,\beta}:=c_\beta^{-1}\, Z_{L,\beta}\,  e^{-\beta L}= 
\sum_{N=1}^{L} \Gamma_\beta^{\, N}\   \mathbf{P}_{\beta}(\cV_{N,L-N}) \quad \text{with}\quad \Gamma_\beta:=\frac{c_\beta}{e^\beta},
\end{align} 
where 
$\mathbf{P}_{\beta}$ is the law of a random walk $V:=(V_i)_{0}^\infty$ starting from the origin ($V_0=0$) and with Laplace symmetric increments, i.e., 
$(V_{i+1}-V_i)_{i\geq 0}$ is an i.i.d. sequence of random variables satisfying 
\begin{equation}\label{lawP}
\mathbf{P}_{\beta}(V_1=k)=\tfrac{e^{-\frac{\beta}{2}|k|}}{c_{\beta}}\quad\forall k\in\mathbb{Z}\quad\text{with}\quad c_{\beta}:=\tfrac{1+e^{-\beta/2}}{1-e^{-\beta/2}}, 
\end{equation} 
and where for every $N\in \{1,\dots,L\}$ the set $\cV_{N,L-N}$ gathers the $N+1$ step trajectories of the random walk  
sweeping a geometric area $L-N$ and finishing at $0$, i.e.,
\be{defG}
\cV_{N,L-N}:=\{V\in \Z^{N+1}\colon G_N(V)=L-N, V_{N+1}=0\} \quad \text{with} \quad  G_N(V)=\textstyle \sum_{i=0}^{N} |V_i|.
\ee
\smallskip
The excess free energy $\tilde f(\beta)$ corresponds to the exponential growth rate of  $\tilde Z_{L,\beta}$ and therefore
can be deduced from the convergence radius of the grand canonical free energy. Using \eqref{tgh} we obtain  
\be{grandcan}
\sum_{L\geq 1} \tilde Z_{L,\beta} e^{-\delta L}=\sum_{N=1}^{\infty} \big(\Gamma_\beta e^{-\delta}\big)^N 
\mathbf{E}_{\beta}\big(e^{-\delta \,  G_N(V)}\,  \ind_{\{V_{N+1}=0\}}\big).
\ee
and since  $\beta\mapsto \Gamma_\beta$ is decreasing on $[0,\infty)$ and satisfies $\Gamma_0>1$ and 
$\lim_{\beta\to \infty} \Gamma_\beta=0$ we deduce from \eqref{grandcan} that the critical point $\beta_c$ of IPDSAW is the unique solution of $\Gamma_\beta=1$. With Theorem \ref{Thm3} below, we derive from  \eqref{grandcan} a simple formulation for the free energy. To that aim we set  
\begin{equation}\label{eq:funch}
h_\beta(\delta):=\lim_{N\to\infty}\frac{1}{N}\log\mathbf{E}_\beta\bigl(e^{-\delta\, G_N(V)}\bigr), \quad \delta \geq 0.
\end{equation}

%From \eqref{tgh}, one derives the phase diagram of the model as follows. By restricting the computation of $Z_{L,\beta}$ to the single trajectory $\ell \in \cL_{\sqrt{L},L}$ (for $L\in \N^2$) which maximizes the Hamiltonian (i.e.,  $\ell_i=(-1)^{i-1} \sqrt{L}-1$ for 
% $i\leq  \sqrt{L}$  we get that $f(\beta)\geq \beta$ for every $\beta\geq 0$.  Since  the sum in \eqref{tgh} is $O(L)$ when $\Gamma_\beta\leq 1$ we deduce that 
%$f(\beta)=\beta$ for $\beta\geq \beta_c$. When $\beta<\beta_c$
%show that when $\gamma_\beta\leq 1$ and therefore $f(\beta)\leq \beta$.  Finally, when $\Gamma_\beta>1$ (i.e., 
%$\beta<\beta_c$) the 
%  denoted by $\beta_c$  One can easily 
%  \begin{equation}\label{eq:main}
%\tilde{f}^\mathsf{m}(\beta)=\sup_{\alpha\in[0,1]}\left[\alpha\log\left(\Gamma^\mathsf{m}(\beta)\right)+\alpha\,g_\beta\left(\tfrac{1-\alpha}{\alpha}\right)\right],
%\end{equation}
%
%\begin{equation}\label{eq:glimit1}
%g_\beta(\alpha):=\lim_{N\to\infty}\tfrac{1}{N}\log\mathbf{P}_{\beta}\big(A_N(V)\leq\alpha N,\,V_N=0\big), \quad  \alpha\in [0,\infty).
%\end{equation}

\begin{theorem}[\cite{CNGP13}, Theorem A]\label{Thm3}
The excess free energy  satisfies
\begin{itemize}
\item   $\tilde{f}(\beta)$ is the unique $\delta$-solution of $\log(\Gamma_\beta)-\delta+h_\beta(\delta)=0$ for  $\beta<\beta_c$
\item  $\tilde{f}(\beta)=0$ for  $\beta\geq \beta_c$
\end{itemize}

%\smallskip

\end{theorem}
Theorem \ref{Thm3} draws a tight link between the asymptotics of $\beta\mapsto \tilde f(\beta)$ at $\beta_c^-$ and the asymptotics of
$\gamma\mapsto h_\beta(\gamma)$ at $0^+$. This is the key to prove Theorem \ref{Thm2} which gives a complete answer to the first challenge raised in Section \ref{challenges} (recall \eqref{asydev}).

\begin{theorem}[\cite{CNGP13}, Theorem B]\label{Thm2}
The collapse transition of IPDSAW   is second order with critical exponent   $3/2$. Moreover, the first order Taylor development of the excess free energy at  $\beta_c^-$ is given by 
\begin{equation}
\lim_{\gep\to 0^+}\frac{\tilde{f}(\beta_c-\epsilon)}{\gep^{3/2}}=\Big(\frac{c}{d}\Big)^{3/2},
\end{equation}
with $\sigma_\beta^2=\mathbf{E}_\beta(V_1^2)$ and 
	$c=1+\tfrac{e^{-\beta_c/2}}{1-e^{-\beta_c}}$,
 and with 
\be{c1}
d=-\lim_{T\to \infty}\frac{1}{T}\log\mathbf{E}\bigr(e^{-\sigma_{\beta_c}\int_0^T|B(t)|dt}\bigr)=2^{-1/3}|a'_1|\sigma_{\beta_c}^{2/3},
\ee
with $a'_1$ the smallest zero (in absolute value) of the first derivative of Airy function.
\end{theorem}

\begin{remark}
The computation of $\mathbf{E}(e^{-s\int_0^1 |B_s| ds})$  for $s>0$ 
is due to \cite{K46} (see e.g. \cite{J07}).
\end{remark}

Let us explain in few words how  Theorem \ref{Thm2} can be deduced from  Theorem  \ref{Thm3}. 
One easily understand that the asymptotic development of   $\tilde{f}(\beta)$ at $(\beta_c)^-$ is strongly related to the fact that there exists a  
constant $c>0$ such that 
\be{expa}
h_\beta(\delta)=-c\,\delta^{2/3}+o(\delta^{2/3}), \quad \text{as $\delta\to 0^+$}.
\ee
These asymptotics are obtained by applying a coarse graining argument: we partition the $V$ random walk trajectory into independent blocks, of size $T\delta^{-2/3}$ with $T\in\mathbb{N}$ chosen arbitrarily and   $\delta$ small enough. Thus, a $N$-step $V$ trajectory is decomposed into  $N/(T\delta^{-2/3})$ blocks that are subsequently used  to prove that as  $\delta\searrow 0$, we have 
\begin{equation}\label{eq:hhbeta}
\lim_{N\to\infty}\frac{1}{N}\log\mathbf{E}_\beta(e^{-\delta G_N})\sim\lim_{T\to\infty}\frac{\delta^{2/3}}{T}\log\mathbf{E}_\beta(e^{-\delta G_{T\delta^{-2/3}}}).
\end{equation}
Donsker's invariance principle ensures  (assuming for simplicity   $\mathbf{E}_\beta(V_1^2)=1$) (cf \cite[p.~405]{RD05}) that 
\begin{equation}\label{ervu}
k^{-3/2}\sum_{i=1}^{Tk}|V_i|\xrightarrow{\text{Law}}\int_0^T|B(t)|dt\quad\text{as}\ k\to\infty,
\end{equation}
where $B$ is a standard Brownian motion. Thus, we choose  $k=\delta^{-2/3}$ in \eqref{ervu} and since  $|e^{-\delta G_{T\delta^{-2/3}}}|\leq1$, we conclude that
\begin{equation}\label{eq:Donsker}
\mathbf{E}_\beta\big(e^{-\delta G_{T\delta^{-2/3}}}\big)\to\mathbf{E}\big(e^{-\sigma_{\beta_c} \int_0^T|B(t)|dt}\big)\quad\text{as}\ \delta\to 0.
\end{equation}
This last convergence combined with \eqref{eq:hhbeta} implies that  $h_\beta(\delta)\sim -c\,\delta^{2/3}$ where $c$ can be expressed using the Laplace transform of the Brownian area, i.e.,
\begin{equation}
c=-\lim_{T\to\infty}\frac{1}{T}\log\mathbf{E}\big(e^{-\sigma_{\beta_c} \int_0^T|B(t)|dt}\big)>0.
\end{equation}

\section{Geometric characterization of IPDSAW}\label{Geocha}

\subsection{Random walk representation}
As mentioned above, the probabilistic expression of the partition function displayed in \eqref{tgh} 
is obtained after mapping appropriately the trajectories of IPDSAW onto random walk trajectories.
 Let us be more specific by recalling \eqref{defLL} and \eqref{defG}  
and, for every $N\in \{1,\dots,L\}$,  by settling a one-to-one correspondance $T_N$ that maps  $\cV_{N,L-N}$ onto $\cL_{N,L}$, i.e., 
\be{defTN}
T_N(V)_i=(-1)^{i-1} V_i\quad  \text{for all} \quad i\in \{1,\dots N\}.
\ee 

%and the polymer measure in \eqref{polmes} becomes
%\begin{equation}
%P_{L,\beta}(l)= \frac{ e^{\beta H_{L}(w)}}{Z_{L,\beta}}, \quad l\in \Omega_L, 
% \end{equation}
%\begin{figure}[ht]\center
%	\includegraphics[width=.3\textwidth]{picstretches.pdf}
%	\caption{Example of a trajectory with $N=5$ vertical stretches
%          and length $L=16$.}
%	\label{fig:stretches}
%\end{figure}
We note that 
$\forall x,y\in\mathbb{Z}$ one can write  $x\;\tilde{\wedge}\;y=\frac12\left(|x|+|y|-|x+y|\right)$ and therefore  the partition function  defined initially in \eqref{defpolme} becomes
\begin{align}\label{ls}
\nonumber Z_{L,\beta}
&=\sum_{N=1}^{L}\sum_{\substack{l\in\mathcal{L}_{N,L}\\l_0=l_{N+1}=0}}\exp{\Bigl(\beta\sum_{n=1}^N{|l_n|}-\tfrac{\beta}{2}\sum_{n=0}^N{|l_n+l_{n+1}|}\Bigr)}\\
&=c_\beta\, e^{\beta L} \sum_{N=1}^{L}\left(\tfrac{c_\beta}{e^\beta}\right)^N\sum_{\substack{l\in\mathcal{L}_{N,L}
\\ l_0=l_{N+1}=0}}\prod_{n=0}^{N}\frac{\exp{\Bigl(-\tfrac{\beta}{2}|l_n+l_{n+1}|\Bigr)}}{c_\beta}.
\end{align}
%At this stage we recall the definition of the auxiliary random walk $V$ in  (\ref{lawP}--\ref{defG}) as well as the family of one to one correspondence $(T_N)_{N=1}^L$
%between path configurations and random walk trajectories (see \ref{defTN}). 
At this stage, we note that  for $l\in \cL_{N,L}$  the increments  of $(V_i)_{i=0}^{N+1}=(T_N)^{-1}(l)$  in \eqref{defTN} satisfy $V_i-V_{i-1}:=(-1)^{i-1}(l_{i-1}+l_i)$. Therefore, 
\begin{align}\label{ls2}
\nonumber \tilde Z_{L,\beta}
&=\sum_{N=1}^{L}\Gamma_\beta^{\,N}\sum_{\substack{l\in\mathcal{L}_{N,L}
\\ l_0=l_{N+1}=0}}\mathbf{P}_\beta(V=(T_N)^{-1} (l))
%&=\sum_{N=1}^{L} \Gamma_\beta^{\,N}\sum_{V\in \cV_{N,L-N}}\prod_{n=0}^{N} \mathbf{P}_\beta(U_{n+1}=V_{n+1}-V_n),
\end{align}
which implies \eqref{tgh}.

Another useful consequence of formula \eqref{tgh} is that it gives a method to sample IPDSAW trajectories 
with the help of random walk paths. To be more specific, let us denote by $N_l$ the horizontal extension of a given $l\in \Omega_L$, i.e.,  $l\in \cL_{N_l,L}$. Since in \eqref{tgh}, the term indexed by $N$  in the summation 
corresponds to the contribution of $\cL_{N,L}$ to the partition function we can state that
\be{lawex}
P_{\beta,L}(N_l=k)=\frac{\Gamma_\beta^k  \mathbf{P}_\beta(\cV_{k,L-k})}{\sum_{k=1}^L \Gamma_\beta^k  \mathbf{P}_\beta(\cV_{k,L-k})}, \quad k\in \{1,\dots,L\},
\ee
and that for every $N\in \{1,\dots,L\}$, 
\begin{align}\label{egalisucc}
%P_{L,\beta}(N_l=N)&={\bf P}_\beta(\xi_L=N\, |\, V_{\xi_L+1}=0, K_{\xi_L}=L), \\
 P_{L,\beta}(l\in \cdot \, |\, N_l=N)&={\bf P}_\beta(T_N(V)\in \cdot\, |\, V\in \cV_{N,L-N}).
\end{align}
As a consequence, one can first sample an extension $N$ under $P_{L,\beta}$ with \eqref{lawex} and then 
sample a $V$ trajectory under  ${\bf P}_\beta$ conditioned on  $\cV_{N,L-N}$ and finally apply $T_N$ to $V$ 
to obtain an IPDSAW trajectory. This method can be implemented to simulate long critical IPDSAW trajectory
(see Fig. \ref{fig:simulations}).

% with a very strong link between the polymer law  $P_{L,\beta}$ and the random walk law ${\bf P}_\beta$
%conditioned on a suitable event. We recall (\ref{taul}--\ref{taul}), the fact that $\Gamma_{\beta_c}=1$ and also that the term indexed by $N$ in the sum in \eqref{ls} corresponds to the contribution to the partition function of those 
%path in $\cL_{N,L-N}$. Consequently, we can derive from (\ref{ls}--\ref{tgh}) that for every $N\in \{1,\dots,L\}$, 
%\begin{align}\label{egalisucc}
%P_{L,\beta}(N_l=N)&={\bf P}_\beta(\xi_L=N\, |\, V_{\xi_L+1}=0, K_{\xi_L}=L), \\
%\nonumber P_{L,\beta}(l\in \cdot \, |\, N_l=N)&={\bf P}_\beta(T_N(V)\in \cdot\, |\, \xi_L=N,  V_{N+1}=0, K_{N}=L).
%\end{align}
%Theorem \ref{transf21} is a straightforward consequence of \eqref{egalisucc}.
%
%
%
%

\subsection{Scaling limit of IPDSAW in each regime}\label{sec:scalingIPDSAW}

To describe geometrically an IPDSAW configuration $l \in \cL_{N,L}\subset \Omega_L$ , one  may  consider its {\it upper envelope} $\cE^+_l$ (respectively  {\it lower envelope} $\cE^-_l$), namely the random process that links consecutively the top (resp. the bottom) of every vertical stretch   
constituting $l$, i.e., 
$\cE^+_{l,0}=\cE^-_{l,0}=0$ and $\cE^+_{l,N+1}=\cE^-_{l,N+1}=l_1+\dots+l_{N}$ and 
\begin{align}\label{trek}
\cE^+_{l,i}&=\max\{l_1+\dots+l_{i-1}, l_1+\dots+l_{i}\},\quad i\in \{1,\dots,N\},\\
\nonumber \cE^-_{l,i}&=\min\{l_1+\dots+l_{i-1}, l_1+\dots+l_{i}\},\quad i\in \{1,\dots,N\}.
\end{align}
%\begin{align}\label{caract}
% \cE^+_{l,i}&=M_{l,i}+\frac{|V_{l,i}|}{2},\quad i\in \{0,\dots,N+1\},\\
%\nonumber \cE^-_{l,i}&=M_{l,i}-\frac{|V_{l,i}|}{2},\quad i\in \{0,\dots,N+1\}.
%\end{align}
Since a given configuration $l$ sampled from $P_{L,\beta}$  fills entirely the subset of $\N\times \Z$ trapped in-between those two envelopes, the scaling limit of IPDSAW (as its length $L$ diverges) is obtained by determining the limiting law of $(\cE_l^-,\cE_l^+)$ rescaled in time and space appropriately.
\smallskip

Another geometric description of $l\in \Omega_L$ can be made  
by considering two auxiliary processes, i.e., the {\it profile}
$|l|:=(|l_i|)_{i=0}^{N+1}$ (with $l_0=l_{N+1}=0$ by convention) and the  center-of mass walk $M_l:=(M_{l,i})_{i=0}^{N+1}$ that links the middle of each stretch consecutively, i.e., $M_{l,0}=0$ and  $M_{l,N+1}=l_1+\dots+l_N$ and 
\be{droitmi}
M_{l,i}=l_1+\dots+l_{i-1}+\frac{l_i}{2},\quad i\in \{1,\dots,N\}.
\ee
Working with $(\cE_l^-, \cE_L^+)$ or with $(|l|,M_l)$ turns out to be equivalent since 
$ \cE_L^+=M_l+\frac{|l|}{2}$ and $ \cE_L^-=M_l-\frac{|l|}{2}$. For simplicity, our results will be displayed with $(|l|,M_l)$
because asymptotically the profile and the center-of-mass always decorrelate.

We define a scaling operator $T_{\alpha,\beta}$ which rescales simultaneously the profile and the center-of-mass walk 
by $L^\alpha$ in time and by $L^\beta$ in space, i.e., 
\begin{equation}
T_{\alpha,\beta}(l)=\frac{1}{L^\beta} \Big(M_{l,\lfloor tL^\alpha\rfloor\wedge N_l},
 |l_{\lfloor tL^\alpha\rfloor\wedge N_l}|\Big)_{t\in [0,\infty ]}.
\end{equation}
Before stating Theorem \ref{theolim} below, we recall that
$\sigma_\beta^2=\mathbf{E}_\beta(V_1^2)$ (see \ref{lawP}). Let us also say that 
in Theorem \ref{theolim} the convergences occur in distribution for cadlag functions on $[0,\infty)$ endowed with the distance of uniform
convergence on every compact subset of $[0,\infty)$. For simplicity, all processes in the statement of Theorem \ref{theolim} have a finite time horizon but we implicitly consider that they remain constant afterwards and therefore are defined on $[0,\infty)$.  Theorem \ref{theolim}
gathers results from \cite[Theorem D]{CNGP13} and \cite[Theorem 2.8]{CP15} and \cite[Theorem C]{CarPet17a}.

\begin{theorem}[]\label{theolim}
For $L\in \N$, we consider an IPDSAW trajectory $l$ sampled from $P_{L,\beta}$. Then,
  \begin{enumerate}[(1)]
  \item if $\beta<\beta_c$, 
\be{cvet}
\lim_{L \to \infty} T_{1,\frac{1}{2}}(l)=\alpha_\beta\,  \big( B_{s},0\big)_{s\in [0,e_\beta)},
\ee
with $e_\beta\in (0,1)$ and $\alpha_\beta>0$ two  explicit constants.  \\
\item if $\beta=\beta_c$, 
\be{cvcrit}
\lim_{L \to \infty} T_{\frac{2}{3},\frac{1}{3}}(l)=\big(D_s, |B_s|\big)_{s\in [0,a_1]},
\ee 
with $B$ and $D$ two independent linear Brownian motions of variance 
$\frac{1}{4} \sigma_\beta^2$ and $\sigma_\beta^2$ respectively, with $a_1$ the time at which the geometric area swept by $B$ reaches $1$, i.e.,  $\int_0^{a_1} \valabs{B_u}\, du = 1$ and with 
$B$ conditioned on the event  
$B_{a_1}=0$.\\

\item  If $\beta>\beta_c$, 
\be{cvco}
\lim_{L \to \infty} T_{\frac{1}{2},\frac{1}{2}}(l)=\big(0,\gamma_\beta(s)\big)_{\, s\in [0,a_\beta]},
\ee
with $a_\beta$ an explicit constant and $\gamma_\beta$ a deterministic Wulff shape given by 
 \be{defgamma}
\gamma_{\beta}(s)=a_\beta \int_0^s L'\big[ (\tfrac12-\tfrac{x}{a_\beta})\,  \tilde h_{0}\big(\tfrac{1}{a_\beta^2},0\big)\big] dx, \quad s\in [0,a_\beta]
\ee
with  
\be{defabeta}
a_\beta=\argmax\big\{a\log\Gamma(\beta)-\tfrac{1}{a}\,\tilde h_0\bigl(\tfrac{1}{a^2},0\bigr)+aL_{\Lambda}\bigl(\tilde H\bigl(\tfrac{1}{a^2},0\bigr)\bigr), \ a\in(0,\infty)\big\},
\ee
where 
\begin{align*}
%L(h):=&\log\mathbf{E}_\beta[e^{h V_1}], \quad h\in(-\tfrac{\beta}{2},\tfrac{\beta}{2})\\
\textstyle L_{\Lambda}(H):=&\int_0^1 \log\mathbf{E}_\beta[e^{(xh_0+h_1) V_1}] dx, \quad H\in \cD
\end{align*}
with $L(h):=\log\mathbf{E}_\beta[e^{h V_1}]$ for $h\in(-\tfrac{\beta}{2},\tfrac{\beta}{2})$ and 
$\mathcal{D}:=\bigl\{H=(h_0,h_1)\colon\, \{h_0,h_0+h_1\}\subset (-\tfrac{\beta}{2},\tfrac{\beta}{2}\bigr)\bigr\}$ and with
$\tilde H=(\tilde h_0, \tilde h_1)$ the inverse function of $\nabla L_\Lambda(H)$ that is a $\cC^1$ diffeomorphism from $\cD$ to $\R^2$.

\end{enumerate}
\end{theorem}
\smallskip

With Theorem \ref{theolim} we observe that the critical regime is characterized by the fact that  the profile and center-of mass walk of a typical IPDSAW configuration display fluctuations of the same order (i.e., $L^{1/3}$). This is indeed not the case  in the extended regime ($\beta<\beta_c$) and inside the collapsed regime ($\beta>\beta_c$) for different reasons.  

When $\beta<\beta_c$ the self-interaction intensity is weak and therefore the qualitative behavior of a typical IPDSAW trajectory is not different from that of the random walk under its uniform measure (i.e., $\beta=0$). To be more specific, the horizontal extension of a typical trajectory is of order $L$ and the vertical stretches are typically of finite size. We will even see in the proof below that the vertical stretches have an exponential tail. As a consequence the profile vanishes when rescaling it in space by any function  growing say faster than $\log L$ whereas the center-of-mass walk  asymptotically decorrelates from the profile and displays Brownian fluctuations.

When $\beta>\beta_c$, in turn, a typical IPDSAW trajectory performs $L (1+o(1))$ self-touchings (saturation) and therefore must be made of few  large vertical stretches with alternating signs. As a consequence the horizontal extension and the vertical stretches of a typical configuration are both of order $\sqrt{L}$. This strong geometric constraint forces the profile rescaled in time and space by $\sqrt{L}$ to converge towards a deterministic Wulff shape (a sketch of the proof is displayed below). The rescaled center-of mass walk vanishes in the limit \eqref{cvco}. The reason is that, the center-of-mass walk
asymptotically decorrelates from the profile and therefore follows the law of a symmetric random walk of length $\sqrt{L}$ with vertical fluctuation 
of order $L^{1/4}$.

\begin{remark}\label{remasy}
\rm Appart from the extended case, the proof of Theorem \ref{theolim} heavily relies on formulas (\ref{lawex}--\ref{egalisucc}) which 
allows us to work with random walk trajectories under a particular conditioning and subsequently to re-express  
the results in terms of IPDSAW via the applications $T_N$ with $N\leq L$ (recall \ref{defTN}). 
\end{remark}
Let us now give the main steps of the proof of Theorem \ref{theolim} in each regime, starting with the collapsed phase.

\subsubsection{Collapsed regime $\beta>\beta_c$.} Rewrite  \eqref{tgh} as 
\begin{align}\label{tgh2}
\tilde Z_{L,\beta}=
\sum_{N=1}^{L} \exp\big(N \big[\log \Gamma_\beta\, + \tfrac{1}{N} \log  \mathbf{P}_{\beta}(\cV_{N,L-N})\big]\big).
\end{align} 
There are two growth rates of $N$ (as a function of $L$) for which $\mathbf{P}_{\beta}(\cV_{N,L-N})$
has a non trivial exponential decay rate (as a function of $N$), namely, $N\sim L$ and $N\sim \sqrt{L}$.
Therefore, and since  $\Gamma_\beta<1$, the sum in \eqref{tgh} is dominated by those terms indexed by $a \sqrt{L}$ with $a\in (0,\infty)$.  Consequently, we set 
\be{def:gbeta}
g_\beta(u):=\lim_{N\to \infty} \mathbf{P}_{\beta} (G_N(V)=u N^2, V_{N}=0),\quad  u\in (0,\infty),
\ee
and $\tilde Z_{L,\beta}$ is well approximated by 
$$\sum_{a\in \frac{\N}{\sqrt{L}}} \exp\big(\sqrt{L} \big[a \log \Gamma_\beta\, + a g_\beta(\tfrac{1}{a^2})\big]\big)$$
so that $a_\beta$ indeed equals 
$\argmax\big\{a \log \Gamma_\beta\, + a g_\beta(\tfrac{1}{a^2}), \, a\in (0,\infty)\big\}$. 

At this stage, proving \eqref{defabeta} 
simply requires  to provide an analytic expression of $g_\beta$.
To that aim, for $u>0$ we observe that  $\{G_N(V)=u N^2, V_{N}=0\}$ is a large deviation event. Its decay rate can indeed be expressed with $J$  the rate function of Mogulskii Theorem  applied to the rescaled process $\tilde V_N:=\big(\tfrac{1}{N} V_{\lfloor sN \rfloor}\big)_{s\in [0,1]}$ viewed as a random element in  $\mathfrak{B}_{[0,1]}$ the set of cadlag real functions on $[0,1]$, endowed with 
the $L^\infty$ norm. Thus $J:  \mathfrak{B}_{[0,1]}\to [0,\infty]$ 
is defined as
\be{defJ}
J(\gamma)=\begin{dcases*}
	\int_0^1 L^*(\gamma'(t)) dt\quad \text{if} \quad \gamma \in \cA\cC,\\
   +\infty \quad \text{otherwise},
  \end{dcases*}
\ee
where $\cA\cC$ is the set of absolutely continuous functions and where $L^*$ is the Legendre transform of $L$. Rewritting 
$\{G_N(V)=u N^2, V_{N}=0\}$ as  $\{ G(\tilde V_N)=u, \tilde V_N(1)=0\}$ (with $G(\gamma)=\int_0^1 |\gamma(s)| ds$)
and applying Mogulskii Theorem in \eqref{def:gbeta} we obtain
\begin{align}\label{gbetaexpr}
g_\beta(u)&=\inf\big\{J(\gamma),\, \gamma \in \mathfrak{B}_{[0,1]},\, G(\gamma)=u, \gamma(1)=0 \big\}
\end{align}
from which we derive the closed formula 
$g_\beta(u)=-u\,\tilde h_0\big(u,0\big)+L_{\Lambda}\big(\tilde H\big(u,0\big)\big)$. The proof of \eqref{defabeta}
is therefore complete and it remains to prove \eqref{defgamma} by observing that  the infimum in \eqref{gbetaexpr} for $u=1/a_\beta^2$ is attained for  $-\gamma^*_{\beta}$ and  $\gamma^*_{\beta}$ defined as 
\be{defgamma}
\gamma^*_{\beta}(s)=\int_0^s L'\big[ (\tfrac12-x)\,  \tilde h_{0}\big(\tfrac{1}{a_\beta^2},0\big)\big] dx, \quad s\in [0,1],
\ee
so that $\gamma_\beta$ simply satisfies $\gamma_\beta(s)=a_\beta \gamma_\beta^*(s/a_\beta)$ for $s\in [0,a_\beta]$.

\subsubsection{Critical regime $\beta=\beta_c$}
The critical regime is the most delicate since the fluctuations of $|l|$ and $M_l$ are of the same magnitudes and must 
therefore be analyzed simultaneously.  

A few more notations are required here. With $V$ a random walk trajectory and with $j,k\in \N$  we associate  $K_j=j+G_j(V)$ (recall \ref{defG})  and $\xi_k:=\inf\{j\geq 1\colon K_j\geq k\}$. We also associate with $V$ an auxiliary process $M:=(M_i)_{i\in \N}$ build with the 
increments of $V$ as follows: $M_0=0$ and for $j\in \N$
\be{defM}
M_j:=\sum_{i=1}^{j-1} (-1)^{i-1} V_i+(-1)^{j-1} \frac{V_j}{2}=\frac{1}{2}\sum_{i=1}^j (-1)^{i-1}\,  (V_i-V_{i-1}).
\ee
Since  $\Gamma_{\beta_c}=1$ the key tool here is the random walk representation (\ref{lawex}--\ref{egalisucc})  which guarantees that 
\be{defcrit}
P_{L,\beta_c}(l\in \cdot)={\bf P}_\beta(T_{\xi_L}(V) \, |\, K_{\xi_L}=L, V_{\xi_L+1}=0).
\ee
A consequence of \eqref{defcrit} is that $T_{\frac{1}{2},\frac{1}{2}}(l)$ with $l$ sampled from $P_{L,\beta_c}$ has the same law as $(|\tilde V|,\tilde M):=(|\tilde V_s|,\tilde M_s)_{s\in [0,\infty)}$ defined as 
\be{defti}
(\tilde V_s,\tilde M_s)=\frac{1}{L^{1/3}} (V_{\lfloor sL^{2/3}\rfloor \wedge \xi_L},M_{\lfloor sL^{2/3}\rfloor \wedge \xi_L} )
\ee 
where $V$ is sampled from ${\bf P}_{\beta_c}$ and conditioned on   $\big\{K_{\xi_L}=L, V_{\xi_L+1}=0\big\}$. Thus, Theorem 
\ref{theolim} (b) can be proven by considering $(|\tilde V|,\tilde M)$.

{\bf Outline of the proof.} The strategy used in \cite{CarPet17a} consists in decomposing every $V$ trajectory into excursions 
$(\mathfrak{E}_r)_{r\in \N}$ away from the origin.  
The fact that the increments  of $V$ follow a symmetric discrete Laplace distribution yields that 
those excursions (in modulus) $(|\mathfrak{E}|_r)_{r\in \N}$ are independent and have the same distribution (except for the very first one). The conditioning 
$\big\{K_{\xi_L}=L, V_{\xi_L+1}=0\big\}$ under which $V$ is considered gives a particular importance to the geometric areas $(X_r)_{r\in \N}$ swept by the excursions. These areas are i.i.d. and heavy tailed random variables so that it suffices to 
consider finitely many excursions  (those sweeping the largest area) to recover a fraction of the path arbitrary close to $1$. For this reason, for $k\in \N$,  we will truncate 
$(|\tilde V|,\tilde M)$ outside the excursions sweeping an area larger than $L/k$.  Since finitely many excursions of $V$ (the largest ones) are required to reconstruct the truncated process
$(|\tilde V_{L,k}|,\tilde M_{L,k})$,  it should be sufficient to prove a convergence in distribution "excursion by excursion" 
to recover the convergence of the whole truncated process. Then, it remains to control the fluctuations of $V$ and $M$ on the small excursions of $V$ in order to check that their contributions to the limiting process vanish as $k\to \infty$.
\smallskip

 Let us be more specific and define the stopping times $(\tau_r)_{r\in\N}$ by the prescription $\tau_0=0$ and
\begin{equation}\label{premm}
  \tau_{r+1} = \inf\ens{i> \tau_r : V_{i-1}\neq 0 \text{ and }
    V_{i-1}V_i \le 0}.
\end{equation}
For every $r\in \N$ we denote by  $|\mathfrak{E}|_r$  the r-th excursion of $V$ in modulus, i.e.,  
\be{defexc}
%\mathfrak{E}_k=(i,V_i)_{i\in \{\tau_{k-1},\dots,\tau_k-1\}} \quad \text{and} \quad
 |\mathfrak{E}|_r=(i,|V_i|)_{i\in \{\tau_{r-1},\dots,\tau_r-1\}},
\ee
and it turns out (see \cite[Proposition 3.1]{CarPet17a}) that provided we transform slightly the law of $V_0$, the sequence 
$(|\mathfrak{E}|_r)_{r\geq 1}$ is i.i.d.  We  introduce for every $r\in \N$ the sum $X_r$ of the length and of the geometric area swept by $|\mathfrak{E}|_r$, i.e., 
%Thus, we define the length and the area swept by the $k$-th excursion (that are  $|\mathfrak{E}|_k$ dependent only) as
\begin{equation}\label{deuxx}
  X_r:= \tau_r -\tau_{r-1} +\big|V_{\tau_{r-1}}\big | + \cdots + \valabs{V_{\tau_r-1}}. 
\ee
With a slight abuse of notation, we will call $X_r$ the geometric area swept by the $r$-th excursion and we define  
$\mathfrak{X}$ a random set of points on $\N_0$ as
\be{defX}
\mathfrak{X}=\{0\}\cup \{X_1+\dots+X_n, n\in \N_0\}.
\ee
For simplicity, we  transform the conditioning under which $(|\tilde V|,\tilde M)$ is considered into $\{L\in \mathfrak{X}\}$. This 
does not change the scaling limit of $(|\tilde V|,\tilde M)$ and lightens the presentation of the proof. Under the conditioning 
$\{L\in \mathfrak{X}\}$ we denote by $v_L$ the number of excursions completed by $V$ when its geometric area reaches $L$. 

\br{remht}
\rm A crucial result at this stage is that $X_1$ is heavy tailed. Deriving a local limit theorem for the geometric area swept by 
a random walk excursion (say with centered increments that have  finite second moments) was an open issue until recently. The reason is that computing the characteristic function of such geometric area is difficult and therefore
Gnedenko's type arguments can not be applied straightforwardly. In \cite[Theorem 1.1]{DKW13}, such a local limit theorem has been derived giving us  $\lim_{n\to \infty} n^{4/3} {\bf P}_{\beta_c}(X_1=n)=c_1$ and  $\lim_{L\to \infty} L^{2/3}{\bf P}_{\beta_c}(L\in \mathfrak{X})=c_2$ with $c_1,c_2>0$. Thus,  by recalling \ref{tgh}, we obtain sharp asymptotics for the critical partition function, i.e., 
$$ Z_{L,\beta_c}=e^{\beta_c L} \frac{c_3}{L^{2/3}}(1+o(1)), \quad \text{with}\ c_3>0 \ \text{explicit}.$$
\er

\noindent {\it Truncation of the profile and center-of-mass walk.} 
As mentioned in the outline, since the  variables $(X_r)_{r\geq 1}$ are heavy tailed, we truncate $(|\tilde V|,\tilde M)$ outside the excursions sweeping an area larger than $L/k$. 
We recall \eqref{defM} and \eqref{premm} and for every $r\in\N$, we let $M^{\text{exc}}(r)$ be the contribution of the $r$-th excursion to the center-of-mass walk, i.e., 
\be{defexcmid}
M^{\text{exc}}(r)=\sum_{i=\tau_{r-1}}^{\tau_r-1} \, (-1)^{i-1} \, V_i.
\ee 
For  $x\in \N$, we truncate $V$ outside the excursions of geometric area larger than $x$ to obtain $(V_x^+(i))_{i\in \N\cup\{0\}}$. Similarly, with the help of \eqref{defexcmid} we define the discrete process $(M^+_{x}(i))_{i\in \N\cup\{0\}}$
% respectively
%$(M^-_{x}(i))_{i\in \N\cup\{0\}}$ 
which remains constant outside 
%(resp. inside)
the excursions of geometric area larger than $x$ and follows the center-of-mass walk elsewhere, i.e., for every $t\in \N$ and $i\in \{\tau_{t-1},\dots,\tau_t-1\}$
\begin{align}\label{tronc3}
M_{x}^+(i)&:=\sum_{r=1}^{t-1}\,  M^{\text{exc}}(r) \, \ind_{\{X_r\geq x\}}+ \Big[\sum_{j=\tau_{t-1}}^{i-1} (-1)^{j-1} V_j + 
(-1)^{i-1} \frac{V_i}{2}\Big] \,  \ind_{\{X_t\geq x\}},\\
\nonumber V_x^+(i)&:=V_i \,  \ind_{\{X_t\geq x\}}. 
%\nonumber M_{x}^-(i)&:=\sum_{r=1}^{t-1} \,  M^{\text{exc}}(r) \, \ind_{\{X_r< x\}} + \Big[\sum_{j=\tau_{t-1}}^{i-1} (-1)^{j-1} V_j + 
%(-1)^{i-1} \frac{V_i}{2}\Big] \,  \ind_{\{X_t< x\}}.
\end{align}
%We observe that $M=M^+_x + M^-_x$. 
%\begin{align}\label{tronc2}
%M_{i,x}^+&:=\sum_{j=0}^{i-1} (-1)^{j-1}\,  V_j\,   \ind_{\{i\in \cY_x\}}+ (-1)^{i-1} \, \frac{V_i}{2} \,  \ind_{\{i\in \cY_x\}},\\
%\nonumber M_{i,x}^-&:=\sum_{j=0}^{i-1} (-1)^{j-1}\,  V_j \,  \ind_{\{i\notin \cY_x\}}+ (-1)^{i-1}\,  \frac{V_i}{2}\,   \ind_{\{i\notin \cY_x\}}.
%\end{align}
%therefore, for every $t\in \N$ and $i\in \{\tau_{t-1},\dots,\tau_t-1\}$ we can rewrite
For $k\in \N$, the  truncated processes $\widetilde V_{L,k}$ and $\widetilde M_{L,k}$ are obtained from  $V^+_{L/k}$ and $M^+_{L/k}$ as in \eqref{defti}.
\smallskip

\noindent {\it Truncation of Brownian motion.} 
We recall the definition of $B$ and $D$ in the statement of Theorem \ref{theolim} (b). As in the discrete case, we truncate  $B$ and $D$ outside the 
excursions of $B$ sweeping a geometric area larger than $1/k$ to obtain $ B^+_{k}$ and $ D^+_{k}$ 
%be the continuous processes that remains constant outside the 
%excursions of $B$ sweeping a geometric area larger than $1/k$ and follows $D$ elsewhere
, i.e.,
\begin{align}\label{defDk}
 D^+_k(s)&=\int_0^{s} \, \ind_{\Gamma_k}(u) \, d D_u,\\
\nonumber  B^+_k(s)&= B_s  \, \ind_{\Gamma_k}(A_s),
\end{align}
where $A_s:=\int_{0}^s |B_s| ds$ is the geometric area swept by $B$ up to time $s$ and where $\Gamma_k :=\ens{u >0: A_{d_u} -A_{g_u} \ge \unsur{k}}$ with
$d_u=d_u(B):=\inf\ens{t>u : B_t=0}$ , $g_u = \sup\ens{t< u : B_t =0}$
so that  $d_u - g_u$ (resp. $A_{d_u} -A_{g_u}$)  is
the length (resp. the geometric area) of the excursion straddling $u$.
%\todo{Fait:Donner une bonne ecriture de $\Gamma_k$}
\smallskip

The proof of Theorem (b) now consists in proving (\ref{fluctdis}--\ref{convtronca}) below. To begin with we must control the fluctuations of $V$ and $M$ outside 
the large excursions of $V$, i.e., prove that for every $\gep>0$ 
\be{fluctdis}
\lim_{k\to \infty} \limsup_{L\to \infty} {\bf P}_{\beta_c}\Big( \sup_{i\leq \xi_L} \big|V_i-V_{L/k}^+(i)\big|+ \big|M_i-M_{L/k}^+(i)\big|\geq \gep L^{1/3}\, \big | \, L\in \mathfrak{X}\Big)=0,
\ee
and similarly for the fluctuations of $B$ and $D$ outside the large excursions of $B$, i.e., 
\be{fluctcont}
\lim_{k\to \infty} \mathbb{P}\Big( \sup_{s\leq a_1} \big|B_s-B_k^+(s)\big|+ \big|D_s-D_{k}^+(s)\big|\geq \gep\Big)=0.
\ee
Then, we must show that for every $k\in \N$  the truncated discrete process $(\widetilde V_{L,k},\widetilde M_{L,k})$
converge in distribution towards its  continuous counterparts $(B_{k}^+,D_k^+)$, i.e., for every $k\in \N$,
\be{convtronca}
\big(\tilde V_{L,k},\tilde M_{L,k}\big)\xrightarrow[L\to \infty]{d} (B_k^+(s),  D_k^+(s))_{s\in [0,a_1]}.
\ee

The proof of (\ref{fluctdis}--\ref{fluctcont}) is displayed in \cite[Sections 5.2 and 5.3]{CarPet17a}. The difficult part consists in controlling the fluctuations of the discrete center-of-mass walk $M$ outside the large excursions of $V$ \big(i.e., of 
$M-M^+_{L/k}$\big). The reason is that at the end of every excursion, i.e. at $\tau_k$ for $k\geq 1$, the 
$V$ trajectory is located very close to the interface, and therefore controlling the fluctuations of $V-V_{L/k}^+$ requires a good  control of the maximum of $V$ on each of its small excursions.  However, this is not the case for the center-of-mass walk, since $M_{\tau_k}$ has no reason to be near the origin. For this reason we must not only control the fluctuations of $M$
inside every small excursions of $V$ but also control the fluctuations of the discrete process of increments 
$\big((M_{\tau_r}-M_{\tau_{r-1}})\,  1_{\{X_r\leq L/k\}}, \ r\leq v_L\big) $.

The proof of \eqref{convtronca} is a reconstruction procedure displayed in \cite[Section 5.1]{CarPet17a}. For every $k \in \N$,  it consists in constructing on the same  probability space a sequence of processes $(Y_L,Z_L)_{L\geq 1}$ and two independent Brownian motions  $B$ and $D$ (with $B$ conditioned on $B_{a_1}=0$) so that 
for every $L\in \N$ the law of $(Y_L,Z_L)$ equals that of $(\tilde V_{L,k},\tilde M_{L,k})$ (with $V$ sampled 
from ${\bf P}_{\beta_c}(\cdot \mid L\in \mathfrak{X}))$ and  $(Y_L,Z_L)$ converges almost surely towards 
$(B_k^+,  D_k^+)$ as $L\to \infty$. To perform this reconstruction we need two auxiliary convergence results. 
First, we recall \eqref{defX} and for every $L\in \N$ we sample $\mathfrak{X}$ from $\mathbf{P}_{\beta}(\cdot\ |\ L\in \mathfrak{X})$. Then,  the tail estimates of $X_1$ (see Remark \ref{remht}) yields that $\frac{1}{L}\mathfrak{X}\cap [0,L]$ converges in law (in the space of closed subsets of $[0,1]$ endowed with the Hausdorff
distance) towards $C_{1/3}\cap[0,1]$ conditioned on $1\in C_{1/3}$ where $C_{1/3}$ is the  $1/3$-stable regenerative set. Second,  we use \cite[Theorem A]{CarPet17b} which yields that, when considering an excursion of the $V$ random walk conditioned on sweeping a prescribed geometric area $L$, the excursion itself and its associated center-of-mass walk, both  rescaled in time by $L^{2/3}$ and in space by $L^{1/3}$  converge in distribution towards $(e_s,D_s)_{s\in[0,a_1]}$ where $e$ is  a Brownian excursion normalized 
by its area and $D$ an independent Brownian motion.

\subsubsection{Extended regime $\beta<\beta_c$}
 The extended regime is somehow the simplest to analyze and was considered in \cite[Section 6]{CP15}.
 % The Random walk representation is indeed 
 %not necessary to provide sharp asymptotics of its partition function and also to derive its scaling limit. 
 The technique consists in partitioning every IPDSAW configuration $l\in \Omega_L$ into $p(l)\in \N$ elementary patterns that do not interact with each other. A pattern is a path whose first zero length vertical stretch occurs only at the end of the path.
 Thus, for $l\in \Omega_L$ we set 
 $ T_0=0$ and  $T_{k}(l) =\inf\big\{j \ge 1+ T_{k-1} : l_j=0\big\}$ for $k\in \{1,\dots,p(l)\}$, 
 so that the total length of the $k$-th pattern in $l$ is  
 \be{defpatt}
 {\textstyle \sigma_k:=T_{k}-T_{k-1}+\sum_{i=T_{k-1}+1}^{T_k} |l_i|.}
 \ee
For $L\in \N$ we denote by $\widehat Z_{L,\beta}$ the contribution to the excess partition function $\tilde Z_{L,\beta}$ (recall \ref{tgh}) of those trajectories that are made of exactly one pattern, i.e., satisfying $\sigma_1(l)=L$. For $\alpha\geq 0$
the moment generating function associated with  $(\widehat Z_{L,\beta})_{L\in \N}$ is 
\be{defphialpha}
\phi(\alpha) := \sum_{L\ge 1} \widehat{Z}_{L,\beta} e^{-\alpha L} \in
]0,+\infty]\,
\ee
and the convergence abcissa
$\hat{f}(\beta) := \inf\ens{\alpha : \phi(\alpha) < +\infty}$. A key observation at this stage is the link between $\phi$ and $\tilde f(\beta)$. It can be proven, using for instance the random walk representation, that 
in the extended phase we have $0 < \hat{f}(\beta) <
\tilde{f}(\beta)$ and moreover that $\phi(\tilde{f}(\beta))=1$. This allows us to define 
the probability $K$ on $\N$ as
\be{defK}
K(n)= \widehat{Z}_{n,\beta}\,  e^{-\tilde{f}(\beta)\,  n}, \quad  n\in \N,
\ee
and to prove that $K$ has an exponential tail.  At this stage, a classical algebraic manipulation of the  partition function $\tilde Z_{L,\beta}$ (see \cite[Section 1.2.1]{cf:Gia}) allows us to show that when $l$ is sampled from $P_{L,\beta}$, the random variables  
$\{\sigma_i(l), i\in \{1,\dots,p(l)\}\}$ are i.i.d. (and conditioned on $\sigma_1+\dots+\sigma_{p(l)}=L$) with law $K$. As a consequence, 
a typical IPDSAW trajectory in the extended phase is made of $O(L)$ patterns of finite size. This yields that 
the vertical stretches of such configurations have length $O(1)$ and therefore explains why the rescaled profile 
in Theorem \ref{theolim} (a) converges to $0$. The Brownian limit of the center-of-mass walk in turn, is easily 
deduced from the fact that, because patterns do not interact energetically, their vertical displacements are 
i.i.d. random variables, centered for obvious symmetry reasons and with a finite second moment 
because the vertical displacement of a pattern is bounded by its total size.

%  and we partition $l$  into $\cup_{k=1}^{\pi_L(l)} \mathfrak{P}_k$ with $\mathfrak{P}_k:=\{l_{T_{k-1}+1},\dots, l_{T_k}\}$
%  and with $\pi_L(l)$ the number of pattern constituting $l$. By construction the Hamiltonian of the path 
%  equals the sum of the Hamiltonians of its patterns and therefore one can compute the grand canonical free energy of the model as 
%  \be{decompozcontraint}
%\tilde{Z}^{\mm,c}_{L,\beta} = \sum_{r=1}^{L/2} \sum_{\substack{s_i\ge 1\\ s_1 +
%  \cdots+s_r=L}} \prod_{i=1}^r \hat{Z}^{\mm,c}_{s_i,\beta}\,.
%\ee
%It is now folklore in Probability theory (see \cite[Chapter 1]{cf:Gia}, for an application of this technique to the linear pinning model) to multiply and divide the r.h.s. by $e^{L \tilde f^\mm(\beta)}$ and obtain
%\be{tre}
%\tilde{Z}^{\mm,c}_{L,\beta} = e^{\tilde f^\mm(\beta)  L}\,  \sum_{r=1}^{L/2} \sum_{\substack{s_i\ge 1\\ s_1 +
%  \cdots+s_r=L}} \prod_{i=1}^r  \PPbeta^{\mm}(\sigma_1=s_i)=\PPbeta^{\mm}(L\in \cT) \, e^{\tilde{f}^\mm(\beta) L}.
%\ee
%  
%
%
%
%

\subsection{More on the collapsed phase: uniqueness of the macroscopic bead.}
Another meaningful manner of describing the geometry of IPDSAW consists in dividing its trajectories into beads. 
More precisely, a bead is made of vertical stretches of strictly positive length and arranged in such a way that two consecutive stretches have opposite directions (north and south) and are separated by one horizontal step (see Fig. \ref{fig:beads}).  A bead ends when the polymer gives the same direction to two consecutive vertical stretches or when a zero length 
stretch appears. 

 Let us define beads rigorously. We pick $l\in \cL_{N,L}$ we set $x_0=0$ and
  $$x_j=\inf\{i\geq x_{j-1}+1\colon\, l_i\;\tilde{\wedge}\;l_{i+1}=0\}\quad \text{for} \quad  j\geq 1.$$
  We let $n(l)$ be the index of the last $x_j$ that is well defined, so that $l$ is partitioned into $\cup_{j=1}^{n(l)} B_j$ where the $j$-th bead $B_j$  is defined as 
\be{defbead}
B_j:=\big\{l_{x_{j-1}+1},\dots,l_{x_j}\big\}, \quad j\leq n(l).
\ee  

Determining the number of beads in a typical configuration is an interesting question raised for instance in \cite{POBG93}. The answer is fairly easy in the extended regime by using the patterns introduced in \eqref{defpatt}. Such a pattern contains at least one bead and since a typical trajectory consists of  $O(L)$ patterns the same 
remains true for the beads.  At criticality, the typical number of beads is $L^{1/3}$. This can be understood 
with the help of the random walk representation (recall (\ref{defcrit})). Every bead of a given $l\in \cL_{N,L}$ is 
indeed associated with an excursion of the associated random walk trajectory $V=(T_N)^{-1}(l)$. Thus the number of beads in a critical IPDSAW trajectory
is also the number of excursions in a $V$ trajectory sweeping a geometric area $L$ and therefore with length $L^{2/3}$ 
(recall Theorem \ref{theolim} (b)). At this stage, the fact that a  $V$ trajectory of length $L^{2/3}$ makes 
order $L^{1/3}$ excursions is sufficient to conclude. In the collapsed regime, the typical number of beads was expected to be small. This has been made rigorous with the following result which states that a typical IPDSAW configuration 
is almost completely trapped inside a {\it unique macroscopic bead}.
\smallskip

For $l\in \cL_{N,L}$ and $j\leq n(l)$ we let $I_j:=|l_{x_{j-1}+1}|+\dots+|l_{x_j}|+x_j-x_{j-1}$ be the size of the $j$-th bead of $l$ (its number of monomers). We also set 
$j_{\text{max}}=\argmax\big\{|I_j|, j\leq n_L(l)\big\}$ so that the size of the largest bead of $l$ is $I_{j_{\text{max}}}$.

\begin{theorem}[\cite{CNGP13}, Theorem C] \label{Thm4}
For  $\beta>\beta_c$, there exists a $c>0$ such that 
\begin{equation}
\lim_{L\to \infty} P_{L,\beta}\big(I_{j_{\text{max}}}\geq L-c\, (\log L)^4\big)=1.
\end{equation}
\end{theorem}

% As proven in \cite{CNGP13} We will prove that the polymer folds itself up into a {\it
%  unique macroscopic bead} and we will identify its horizontal
%extension and its asymptotic deterministic shape.  To quantify these results we need the following notations.
%
\begin{figure}[ht]\center
	\includegraphics[width=.55\textwidth]{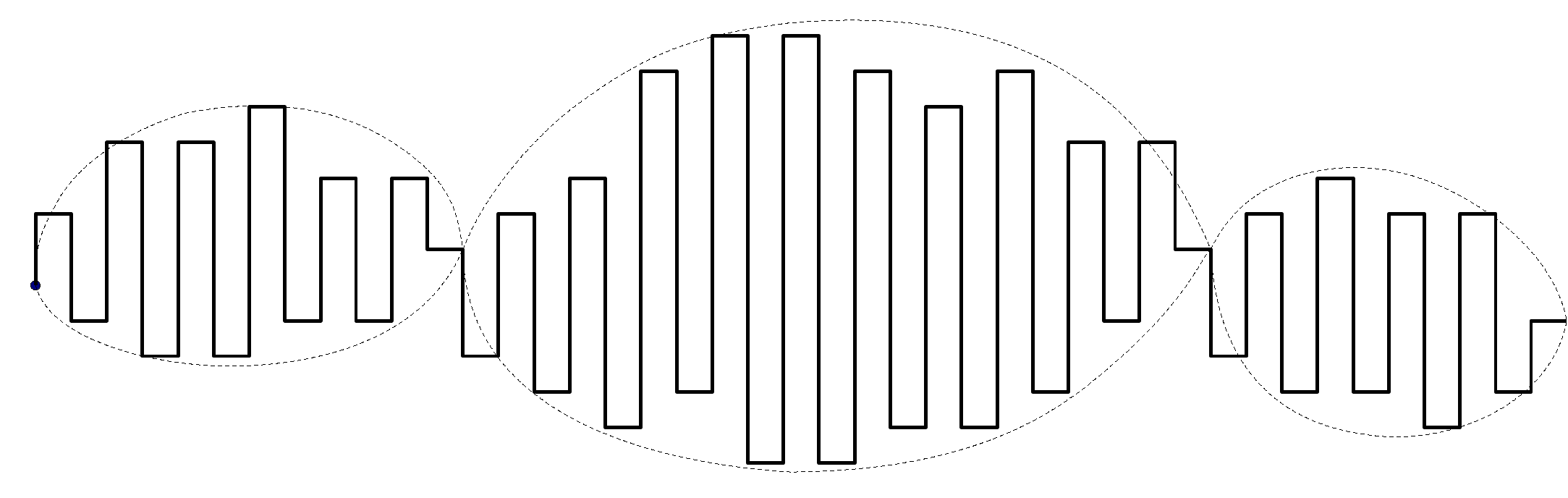}
	\caption{Example of a trajectory with $3$ beads.}
	\label{fig:beads}
\end{figure}

%\subsubsection{Number of beads}
%Let $l\in \Omega_L$  and denote by  $N_L(l)$ its horizontal extension, i.e., $N_L(l)$ is the integer $N$ such that $l\in \cL_{N,L}$.
%Pick $l\in \cL_{N,L}$ and let $(u_j)_{j=1}^N$ be the sequence of cumulated lengths of the polymer after each vertical stretch, i.e.,  $u_j=|l_1|+\dots+ |l_j|+j$ for $j\in \{1,\dots,N\}$. For convenience only, set $l_{N+1}=0$. Set also $x_0=0$ and for $j\in \N$ such that $x_{j-1}<N$, set  (see Fig. \ref{fig:transform}). Finally, let  $n_L(l)$ , i.e., $x_{n_L(l)}=N$. Thus we can decompose any trajectory $l\in \Omega_L$ into a succession of $n_L(l)$ beads, each of them being associated with a subinterval of  $\{1,\dots,L\}$ written as 
%\begin{equation}\label{ut}
%I_j=\{u_{x_{j-1}}+1,\dots,u_{x_j}\}, \quad\text{for}\quad j\in \{1,\dots,n_L(l)\},
%\end{equation}
%and therefore, we can partition $\{1,\dots,L\}$ into $\cup_{j=1}^{n_L(l)} I_j$.
%At this stage, we can define the largest bead of a trajectory $l\in \Omega_L$ as $I_{j_{\text{max}}}$ with
%\begin{equation}\label{ut2}
%j_{\text{max}}=\argmax\big\{|I_j|, j\in \{1,\dots,n_L(l)\}\big\}.
%\end{equation}
%With Theorem \ref{Thm4} below, we claim that, in the collapsed phase, there is only one macroscopic bead. 
%

Theorem \ref{Thm4} can be deduced from Propositions \ref{prop1} and \ref{lem2} below. Proposition \ref{prop1} gives bounds on the partition function restricted to those trajectories describing one bead only. We will discuss the proof of Proposition \ref{prop1} at the end of the present section. Proposition \ref{lem2} shows that the horizontal extension of a typical configuration inside the collapse phase is bounded above by $a\sqrt{L}$ for some $a>0$. This has been discussed already in (\ref{tgh2}-- \ref{def:gbeta}).

% be deduced from 	  Then, we 
%will discuss the proof 
%
% Then we will discuss rye 
%
%Proving Theorem \eqref{Thm4} requires compute the contribution to the partition function of those $l\in\Omega_L$ that have only one bead, i.e. $n_L(l)=1$.  Then, using the fact the Hamiltonian of a configuration 
%

%   to split the partition function of IPDSAW depending 

The subset of $\Omega_L$ containing the  one bead trajectories  is $\Omega^\mathsf{o}_L:=\bigcup_{N=1}^L\mathcal{L}^\mathsf{o}_{N,L}$, with %where $\mathcal{L}^\mathsf{o}_{N,L}$ is the subset of $\mathcal{L}_{N,L}$ defined as
\begin{equation}\label{deftu}
\textstyle\mathcal{L}^\mathsf{o}_{N,L}=\bigl\{l\in\mathcal{L}_{N,L}\colon\,l_i\;\tilde{\wedge}\;l_{i+1}\neq 0\ \forall i\in\{1,\dots,N-1\}\bigr\},
\end{equation}
and its contribution to the partition function is $Z^\mathsf{o}_{L,\beta}:=\sum_{l\in\Omega_L^\mathsf{o}}e^{\beta H_{L}(l)}$.
\smallskip

\begin{proposition}[\cite{CNGP13}, Proposition 4.2]\label{prop1}
For  $\beta>\beta_c$, there exist $c,c_1,c_2>0$ and $\kappa>1/2$ such that 
\begin{equation}\label{state}
\frac{c_1}{L^\kappa}\,e^{\beta L-c\sqrt{L}}\leq Z_{L,\beta}^{\mathsf{o}}\leq\frac{c_2}{\sqrt{L}}\,e^{\beta L-c\sqrt{L}},\quad L\in \N.
\end{equation}
\end{proposition}
\smallskip

\begin{proposition}[\cite{CNGP13}, Lemma 4.1]\label{lem2}
For  $\beta>\beta_c$, there exist $a,a_1,a_2>0$ such that
\begin{equation}\label{rst}
P_{L,\beta}(N_L(l)\geq a_1\sqrt{L})\leq a_2\,e^{-a \sqrt L},\quad L\in \N.
\end{equation}
\end{proposition}

Propositions \ref{prop1} and \ref{lem2} are sufficient to prove Theorem \ref{Thm4}. 
The first step consists in showing that there exists an $s>0$ such that  a typical configuration has exactly one  bead larger than $s (\log L)^2$. To that aim, we use Proposition \ref{prop1} combined with the inequality 
$\sqrt{x}+\sqrt{y}-\sqrt{x+y}\geq \frac{1}{2} \sqrt{x\wedge y}$ to assert  that a typical trajectory has at most one  bead larger than $s (\log L)^2$. 
Then, we note that each bead contains at least one horizontal step. Therefore 
a configuration that has no bead larger than $s (\log L)^2$ has at least  $L/s (\log L)^2$
horizontal steps and can not be typical because it would contradict Proposition \ref{lem2}. 

The second step consists in bounding above the number of monomers that do not belong to the unique big bead. 
%to prove that there exists a $c>0$ such that there are typically  less than $c (\log L)^4$ monomers outside this unique big bead. 
We denote by  $x_1$ (resp. $x_2$) the number of monomers  before (resp. after) this big bead and we assume for instance that $x_1\geq c (\log L)^4$. There are at least $x_1/ s (\log L)^2$ beads between $0$ and $x_1$ and consequently at least 
$x_1/s (\log L)^2$  horizontal steps. The fact that $x_1$ is the end of a bead allows us to split the path at $x_1$ 
and to focus on the trajectory between $0$ and  $x_1$. Applying Proposition \ref{lem2} with $L=x_1$ yields that, 
 typically, less than  $a_1\sqrt{x_1}$ of the first $x_1$ steps are horizontal. This provides a contradiction because, 
 by choosing $c$ large enough, the inequality  $x_1\geq c (\log L)^4$ yields that $a_1\sqrt{x_1}=o(x_1/s (\log L)^2)$.
The proof of Theorem \ref{Thm4} is complete.
 
 \smallskip

Let us give some insights concerning the proof of Proposition \ref{prop1}. The key tool here is the link between 
beads of IPDSAW and excursions of the $V$ random walk. To be more specific, for $k \leq n$,  we 
let  $\mathcal{V}_{n,k}^{+}$ be the subset containing those positive excursions of the $V$ random walk, returning to the origin after $n$ steps, and sweeping an area $k$, i.e., 
\begin{equation}\label{defVV+}
\mathcal{V}_{n,k}^+:=\{V\colon\,V_n=0,\,G_n=k,\,V_i>0\ \forall i\in\{1,\dots,n-1\}\}.
\end{equation}
By mimicking \eqref{tgh} and by noticing that the $T_N$-transformation is a one-to-one correspondence between 
$\mathcal{V}^+_{N+1,L-N}$ and  $\mathcal{L}^\mathsf{o}_{N,L}$ we obtain  that
\begin{equation}\label{eq:dist11}
\widetilde Z^\mathsf{o}_{L,\beta}:= \frac{1}{2  c_\beta}\, e^{-\beta L} Z^\mathsf{o}_{L,\beta}=\sum_{N=1}^{L} \Gamma_\beta^N\mathbf{P}_{\beta}(\mathcal{V}^+_{N+1,L-N}).
\end{equation}
As discussed in (\ref{tgh2}--\ref{def:gbeta}) the sum in \eqref{eq:dist11} is dominated by those terms 
indexed by $N\sim a\sqrt{L}$. Therefore, proving Proposition \ref{prop1} requires the derivation of some sharp bounds on 
$\mathbf{P}_{\beta}(\mathcal{V}^+_{n, a n^2})$ for $a>0$ and $n\in \N$. By using  tilting 
techniques from \cite{DH96} one obtains local limit theorems for $\mathbf{P}_{\beta}(A_n=a n^2, V_n=0)$
that is  with a bridge instead of an excursion and with the algebraic area $A_n(V):=\sum_{i=1}^n V_i$
instead of the geometric area. This provides the upper bound in Proposition \ref{prop1}. To derive the lower bound in 
Proposition \ref{prop1}, it remains to bound from below the probability that a $V$ trajectory remains positive when conditioned on $A_n=a n^2, V_n=0$ and this is the object of \cite[Proposition 2.5]{CNGP13}.

\section{Open problems}\label{open}
 Open issues related to IPDSAW are numerous. Without pretending to be exhaustive, let us 
display four research directions that are both relevant from a physical standpoint and challenging mathematically. 
\begin{enumerate}
\item {\bf Disordered IPDSAW.}  Taking into account an inhomogeneous solvent and/or 
a copolymer (instead of an homopolymer) when studying the collapse transition phenomenon would be 
an important improvement. With a copolymer, this could be achieved by introducing a random component in the self-touching intensity 
involving monomers $i$ and $j$. One could replace $\beta$ by $\beta+s\xi_{i,j}$ 
with $\{\xi_{i,j}, (i,j)\in \N^2\}$ an i.i.d. field of random variables and $s>0$ a tuning parameter.
In this framework, it would be particularly interesting to investigate the relevance of disorder $\xi$ (Harris criterion), 
that is to figure out whether an arbitrary small $s>0$ rounds the phase transition or not.

\smallskip

\item {\bf Towards 2-dimensional ISAW.}  
One could consider an enhanced version of IPDSAW, i.e., a 2-dimensional model 
whose allowed configurations are not directed anymore. Of course the ISAW itself would be an ideal
choice but we have seen that it is very hard to analyse at a rigorous level. Thus, 
one could consider a model that interpolates between IPDSAW and ISAW,
 in the sense that its allowed configurations are not directed anymore and have a connective constant 
 strictly between that of IPDSAW trajectories and that of ISAW trajectories. 
This is the case for the Interacting Prudent Self-Avoiding Walk that has been investigated recently in \cite{PT17}
and will be discussed further in Section \ref{secIPRSAW} below.
\smallskip

\item {\bf Higher dimensions.} 
Obtaining rigorous mathematical results about 
the collapse transition of a three-dimensional extension of IPDSAW  would be of great interest.
One option would be to consider the above mentioned Interacting Prudent Self-Avoiding Walk in 
dimension 3.

\smallskip

\item {\bf Collapse with adsorption.}
The repulsion of monomers exerted by a poor solvent  and 
the adsorption of those monomers along a hard wall are among the most basic 
interactions between an homopolymer and the medium around it (see \cite{Flory}).
Therefore, building a mathematical model taking both interactions into account is physically appealing. 
This has been done in dimension 2 for instance in \cite{F90} or \cite{FY91}  where the IPDSAW is perturbed by 
a pinning interaction at the $x$-axis that plays the role of an impenetrable horizontal interface.
There are precise conjectures concerning the phase diagram of this model (see \cite{FY91}). 
It is expected to be partitioned into 3 phases (collapsed, extended and pinned) separated by 3 critical curves and 
meeting at one tricritical point. So far, the boundary of the collapsed phase has been computed (see \cite{F90}) but
the rest of the phase diagram lacks rigorous mathematical proofs. 
\end{enumerate}

\section{A non directed model of ISAW:  the IPSAW}\label{secIPRSAW}

In the spirit of the second class of open problems mentioned  in Section \ref{open}, the Interacting Prudent Self-Avoiding Walk (IPSAW) is studied in \cite{PT17}.  
%a model of interacting {\emph non-directed}
%self-avoiding walk. 
In size $L\in \N$, the set of configurations consists of the $L$-step prudent paths
introduced in \cite{TD87b}, i.e., 
\begin{align}
\nonumber \Omega_L^{\tx{PSAW}}=\big\{w:=(w_i)_{i=0}^L\in (\mathbb{Z}^2)^{L+1}\colon\,& w_{0}=0,\  w_{i+1}-w_i\in \{\leftarrow,\rightarrow,\downarrow,\uparrow\},\  0\leq i \leq L-1,\\
&w \  \text{satisfies the prudent condition}\big\},
\end{align}
where the \emph{prudent condition} for a path $w$ means that it does not take any step in the direction 
of a lattice site already visited (see Figure \ref{fig:PWandNEPW}).  We define also $\Omega_L^{\tx{NE}}$ as the subset of $\Omega_L^{\tx{PSAW}}$ containing 
those trajectories with a general north-east orientation (see Figure \ref{fig:PWandNEPW}),
that is, all the prudent trajectories that do not take any step in the direction 
of the set $(-\infty,0]^2$.
 We observe that partially-directed self-avoiding paths
(recall \ref{defLL}) are 
in particular north-east prudent paths and that prudent paths are in particular self-avoiding paths, so that 
\be{include}
\Omega_L\subset \Omega_L^{\tx{NE}}\subset \Omega_L^{\tx{PSAW}}\subset  \Omega_L^{\tx{SAW}} , \quad L\in \N,
\ee
where $\Omega_L^{\tx{SAW}}$ denotes the set of $L$-step self-avoiding paths in $\Z^2$ taking unitary steps.

\begin{figure}[ht]
\centering
\begin{subfigure}[t]{0.4\textwidth}
{\includegraphics[scale=1]{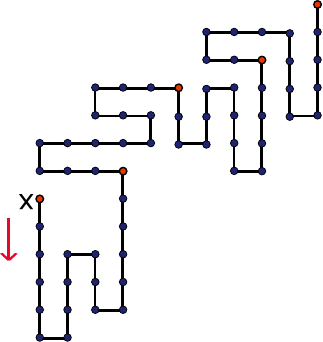}}
\caption{NE-IPSAW}\label{figb}
\end{subfigure}
 \begin{subfigure}[t]{0.4\textwidth}
{\includegraphics[scale=1.0]{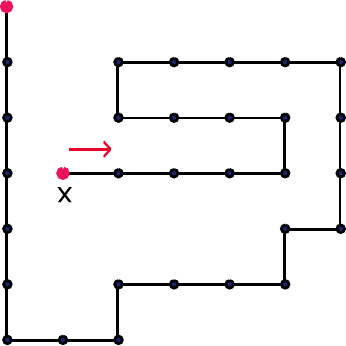}}
\caption{IPSAW}\label{figc}
\end{subfigure}
\caption{
Examples of NE-PSAW ($\text a$) and PSAW ($\text b$) path. 
Each path starts at $x$ and its orientation is 
given by the arrow.
%In ($\text A$) we have drawn an IPDSAW path  
%made of $11$ stretches: $\ell_1=9,\, \ell_2=-7,\, \ell_3=9
%,\, \ell_4=0,\, \ell_5=-12,\, \ell_6=0,\, \ell_7 = 5,\, \ell_8=0,\, \ell_9=5, \,
%\ell_{10}=-7, \ell_{11}=0$. That path performs 
% $19$ self-touching (drawn in red).
}
\label{fig:PWandNEPW}
\end{figure}

To help define the Hamiltonian, 
we associate with every path $w\in \Omega_L^{\tx{SAW}}$ the sequence of those points in the middle of each step, i.e., $u_i=w_{i-1}+\frac{w_i-w_{i-1}}{2}$ ($1\leq i\leq L$). The self-touchings performed by $w$ correspond to the 
non-consecutive pairs $(u_i,u_j)$ at distance one, i.e, $\|u_i-u_j\| = 1$, see Figure \ref{fig:IPSAW}.  Then, the Hamiltonian of every  $w\in \Omega_L^{\tx{SAW}}$ is defined as
\begin{equation}\label{eq:ham1}
\mathrm H_L(w):=\sum_{\substack{i,j=0\\i<j}}^L\ind_{\{\lVert u_i-u_j\rVert=1\}}.
\end{equation}
The coupling parameter $\beta\geq 0$ stands for the self-interaction intensity and 
therefore, the partition functions of IPSAW  and of the North-East model are respectively given by
\begin{equation}
\mathrm{Z}^{\tx{IPSAW}}_{\beta,\, L}:=\sum_{w\in \Omega_L^{\tx{PSAW}}} e^{\, \beta \,   \mathrm H_L(w)} \quad \text{and}  \quad \mathrm{Z}^{\tx{NE}}_{\beta,\, L}:=\sum_{w\in \Omega_L^{\tx{NE}}} e^{\, \beta \,   \mathrm H_L(w)},
\end{equation}
and there exponential growth rate (free energies) by
% by  We denote by $\mathrm F^{\tx{IPSAW}}(\beta)$ and $\mathrm F^{\tx{NE}}(\beta)$ the free energies of the Prudent and Noth-East models, i.e.,
\begin{equation}\label{eq:FreeEnergyIPSAW}
\mathrm F^{\tx{IPSAW}}(\beta):= \lim_{L\to\infty} \frac{1}{L}\log \mathrm Z_{\beta,\, L} \quad \text{and} \quad \mathrm F^{\tx{NE}}(\beta):= \lim_{L\to\infty}\frac{1}{L}\log \mathrm Z^{\tx{NE}}_{\beta,\, L}.
\end{equation}

\begin{figure}[ht]
\centering
\includegraphics[scale=1]{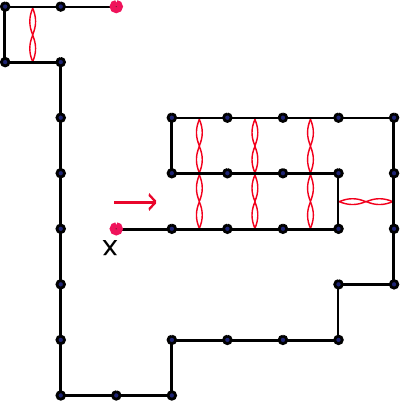}
\caption{A prudent path $w$ in which we highlight the self-touching.
The path starts at $x$ and its orientation is 
given by the arrow. 
In this example we have $L=34$ and 
$\mathrm H_L(w)=8$.
}
\label{fig:IPSAW}
\end{figure}

\br{rem:whyweuse}
\rm Defining the self-touchings of an IPDSAW configurations with pairs of edges or with 
pairs of sites is equivalent (see Fig. \ref{fig:stretches}). This is not the case for non-directed models. The choice "by edge" made in \eqref{eq:ham1} 
for IPSAW turns out to be more tractable than the choice "by site", but
% obtained by defying the Let us stress that the Hamiltonian $H_L(w)$ in  gives rise to the same IPDSAW type of model introduced in Section \ref{S1.1}. 
%In the non-oriented case the choice made in \eqref{eq:ham1} turns out to be more tractable than the choice in which 
%we allow the self-touching between pair of points.
we expect that both choices lead to similar qualitative behaviors.
\er

\br{rem:prudentscal}
\rm The connective constant of North-East prudent paths  $\mathrm F^{\tx{NE}}(0)$ was computed in 
 \cite{B10} (under the name {\it two-sided})  and turns out to be strictly larger than that of partially-directed self-avoiding paths $f(0)$ (recall \ref{defpolme}). This gives an incentive for studying IPSAW since it takes into account much more trajectories than IPDSAW and therefore, the proof of the existence of the collapse transition for IPSAW (see Theorem \ref{Thm1}) is a real step forward. However, one should also acknowledge that 
 IPSAW is still far from ISAW itself. It appears clearly with Theorems \ref{Thm1} and \ref{Thm33} since the collapse transitions of IPSAW 
 and of ISAW are of different nature (see the discussion before Theorem \ref{Thm33} below).   
  \er  

\br{rem:prudent}
\rm  The 
 scaling limit of the 2-dimensional prudent walk has been derived  in  \cite{BFV10} for its kinetic version
 and in \cite{PTS17} with the uniform law. The prudent walk has also been used in \cite{BI15} 
 to build and investigate a non-directed model of polymer adsorption. 
 %the IPSAW for which the  set of allowed spatial configurations for the polymer is given by $\Omega_L^{\tx{PSAW}}$ and its North-East counterpart (NE-IPSAW)
%for which the  set of  configurations  is given by $\Omega_L^{\tx{NE}}$. For both models, each step of the walk is an abstract monomer and 
% we want to take into account the repulsion between monomers and the environment around them. 
% %Physically if such repulsion is strong enough, then the polymer collapses in a small region of the space, 
% %minimizing the contact with the solvent. 
%This is achieved indirectly, by encouraging monomers to attract each other,
%  i.e.,  by assigning an energetic reward $\beta\geq 0$  to any pair of non-consecutive 
%  steps of the walk though adjacent on the lattice $\mathbb{Z}^2$. 
\er

%and the polymer measure takes the form 
%\begin{equation}
%\label{eq:RM}
%\mathrm P_{\beta,\, L}(w)=\frac{e^{\, \beta\,\mathrm H\,(w)}}{\mathrm Z_{\beta,\, L}},\qquad w\in \Omega_L,
%\end{equation}
%where $\mathrm{Z}_{\beta,\, L}$ is the partition function.
%The probability measure \eqref{eq:RM} defines the IPSAW model. 

\smallskip

%The key objects of our analysis are the free energies of both models, i.e., $\mathrm F(\beta)$ and  $\mathrm F^{\tx{NE}}(\beta)$ which record the exponential growth rate of the partition function sequences $(\mathrm{Z}_{\beta,\, L})_{L\in \N}$ and $(\mathrm{Z}^{\tx{NE}}_{\beta,\, L})_{L\in \N}$, respectively.  Thus, 
%\begin{equation}
%\label{eq:FreeEnergyIPSAW}
%\mathrm F(\beta):= \lim_{L\to\infty}\frac{1}{L}\log \mathrm Z_{\beta,\, L} \quad \text{and} \quad \mathrm F^{\tx{NE}}(\beta):= \lim_{L\to\infty}\frac{1}{L}\log \mathrm Z^{\tx{NE}}_{\beta,\, L}.
%\end{equation}
%The convergence in the right hand side of \eqref{eq:FreeEnergyIPSAW} will be proven in Section \ref{app:freeenergy}. The convergence in the l.h.s. of  \eqref{eq:FreeEnergyIPSAW} is more complicated and  it will   be 
%obtained as a by-product of Theorem \ref{thh1} below.

\subsection{Existence of a collapse transition of IPSAW}
\label{MainRes}
 
The prudent model is, to our knowledge, the only 
\emph{non-directed} model of a 2-dimensional interacting self-avoiding walk 
for which the existence of a collapse transition has been proven  rigorously. This is the main result in
\cite{PT17} along with the equality between both free energies in \eqref{eq:FreeEnergyIPSAW} which (at $\beta=0$) answers an open
question raised in \cite{B10}.  

 %This procedure shows that the loss of energy between a prudent path and a $2$-sided prudent path is sub-exponential
%and that the number of prudent paths that can be mapped on a given NE-prudent path by $M_L$ is sub-exponential. 
%At $\beta=0$, for showing that the exponential growth rate of the set of $2$-sided prudent paths
%equals that of generic prudent path we only need to count the number of ancestors by $M_L$. 
%We indeed build a sequence of path transformations $(M_L)_{L\in \N}$ such that 
%for every $L\in \N$, $M_L$ maps any generic path in $\Omega_L^{\tx{PSAW}}$ onto a \emph{2-sided} prudent path in $\Omega_L^{\tx{NE}}$ and  satisfies the following properties:
%\begin{itemize}
%\item for every $w\in \Omega_L^{\tx{PSAW}}$,  the difference between the Hamiltonians of $w$ and of $M_L(w)$ is $o(L)$,
%\item the number of ancestors of a given path in $\Omega_L^{\tx{NE}}$ by $M_L$ can be shown to be $e^{o(L)}$.
%\end{itemize} 
%
%Such a mapping allows us to prove the following theorem. 
\begin{theorem}[\cite{PT17}, Theorem 2.1] \label{thh1} For $\beta\geq 0$, 
\begin{equation}\label{thh}
\mathrm  F^{\rm\tx{IPSAW}}(\beta)= \mathrm  F^{\rm\tx{NE}}(\beta).
\end{equation}
\end{theorem}
%The free energy equality in \eqref{thh} will subsequently be used to establish Theorem \ref{Thm1} below, which states that 
%IPSAW undergoes a collapse transition at finite temperature. 
%Moreover Theorem \ref{thh1} 
%at $\beta=0$ answers an open question from combinatorics. The \emph{2-sided} prudent trajectories have indeed been already studied in 
%the mathematical litterature, see e.g., Bousquet Melou \cite{B10},  Detheridge and Guttman \cite{DG08} or  Beaton and Iliev \cite{BI15}. It was conjectured in  \cite{B10} or \cite{DG08} that the exponential growth rate of the set of $2$-sided prudent paths
%equals that of generic prudent path and this is precisely what Theorem \ref{thh1} says at $\beta=0$. 

% will subsequently be used to establish Theorem \ref{Thm1} below, which states that IPSAW undergoes a collapse transition at finite temperature. 

\smallskip

%\noindent \emph{Collapse transition.}  
%Proving that a $2$-dimensional model of an interacting self-avoiding walk displays a collapse transition 
%when the coupling parameter reaches a critical threshold is known to be difficult.
% We will discuss more in depth this issue in Section \ref{bar} below, but let us already say
%that a rigorous proof of such transition was so far  available for directed models only,  that is with random walks taking steps among three directions (e.g. in $\{\uparrow,\to,\downarrow\}$  
%for the partially directed self-avoiding walk). 
%For this reason, and up to our knowledge, Theorem \ref{Thm1} below is the first proof of the existence of a collapse transition for a non directed model of interacting self-avoiding walk in 
%dimension $2$. 
  \begin{theorem}[\cite{PT17}, Theorem 2.2]
\label{Thm1}
There exists  a $\beta_c^{\tx{IPSAW}}\in (0,\infty)$ such that 
\begin{align}
\mathrm  F^{\rm\tx{IPSAW}}(\beta)&>\beta \quad  \text{for every} \quad  \beta<\beta_c^{\tx{IPSAW}},\\
\nonumber\mathrm   F^{\rm\tx{IPSAW}}(\beta)&=\beta \quad  \text{for every} \quad  \beta\geq \beta_c^{\tx{IPSAW}}.
\end{align} 
Thus, the phase diagram $[0,\infty)$ is partitioned into 
a \emph{collapsed phase}, $\mathcal{C}:= [\beta_c^{\tx{IPSAW}},\infty)$ inside which the free energy \eqref{eq:FreeEnergyIPSAW} 
 is linear and an \emph{extended phase}, $\mathcal{E}=[0,\beta_c^{\tx{IPSAW}})$. 
\end{theorem}
The proof of Theorem \ref{thh1} is purely combinatorial.  It consists  in 
building a sequence of path transformations $(M_L)_{L\in \N}$ such that 
for every $L\in \N$, $M_L$ maps  $\Omega_L^{\tx{PSAW}}$ onto $\Omega_L^{\tx{NE}}$ and  satisfies the following properties:
\begin{itemize}
\item for every $w\in \Omega_L^{\tx{PSAW}}$,  the difference between the Hamiltonians of $w$ and of $M_L(w)$ is $o(L)$,
\item the number of ancestors of a given path in $\Omega_L^{\tx{NE}}$ by $M_L$ can be shown to be $e^{o(L)}$.
\end{itemize} 
The prudent condition guaranties that every  $w\in \Omega_L^{\tx{PSAW}}$ can be decomposed in a unique manner into 
at most $\sqrt{L}$  two-sided subpaths that are either North-East, South-East,  North-West or South-West. The procedure encoded in $M_L$ 
consists in detaching one by one the two-sided blocks composing $w$. This can be achieved by loosing at most $o(L)$ self-touchings. Then, using rotations and symmetries we concatenate the two-sided blocks so as to recover an $L$-step North-East path.

\begin{figure}
\includegraphics[scale=1.1]{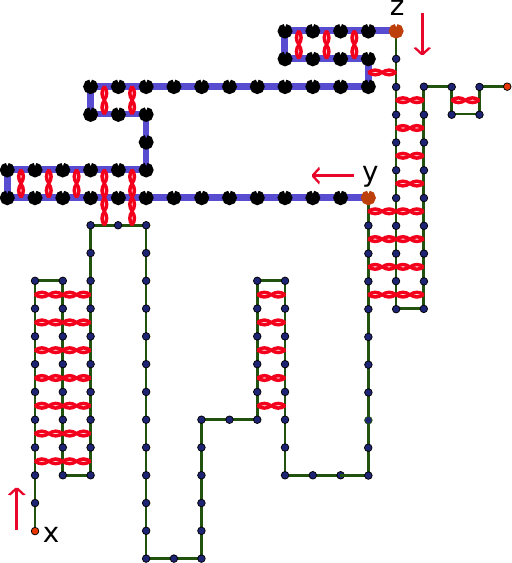}
\caption{A NE-prudent path made of three blocks. 
The first block starts at $x$, the second block starts at $y$ and the third block starts at $z$. Their orientation 
is given by the arrow.
Interactions in each block and between different blocks are highlighted in red. 
%We observe that the $i$-th oriented block of a North-East path may only interact with the  $(i-2)$-th  and $(i-1)$-th blocks. 
%Since we count only one time the self-touching between different oriented blocks, we have that the first stretch of the second block interacts with the inter-stretches of the 
%first block.
We observe that, the first stretch of the third block interacts with both the inter-stretches of the 
second block and the last stretch of the first block.
}
\label{fig:IPDRW2}
\end{figure}

Thanks to Theorem \ref{thh1}, it is sufficient to prove Theorem \ref{Thm1} with the free energy of the North-East model.
The prudent condition yields that every North-East path can be decomposed in a unique way into partially directed subpath (referred to as oriented blocks). Therefore, deriving a sharp upper bound on the North-East partition function $Z_{L,\beta}^{\tx{NE}}$ requires to 
bound from above 
\begin{itemize}
\item the free energy of an oriented block of a given length,
\item the self-touchings between different oriented blocks, 
\item the entropy carried by the fact that the number of oriented blocks and
their respective lengths may fluctuate.
\end{itemize}
Controlling the free energy of an oriented block is achieved with the random walk representation since oriented blocks are 
IPDSAW trajectories. Dealing with the self-touchings occurring between blocks requires to observe that the $i$-th oriented block of a North-East path may only interact with the  $(i-2)$-th  and $(i-1)$-th blocks. Moreover, self-touchings
may only appear between the first stretch of the $i$-th block and the inter-stretches of the 
$(i-1)$-th block or between the first stretch of the $i$-th block and the last stretch of the $(i-2)$-th block (see Figure \ref{fig:IPDRW2}). This allows us to derive an explicit upper bound  on the total number of self-touchings that may appear between the different oriented blocks of a North-East configuration, see Figure \ref{fig:IPDRW3}. Therefore, it remains to
control the entropy related to the fact that the lengths of oriented blocks may fluctuate. This is again taken care of with the random walk representation. To be more precise, at the end of the proof we derive a random walk representation for the whole North-East path and this is sufficient to conclude that for $\beta$ large enough, the free energy is not larger than $\beta$.

\begin{figure}
\includegraphics[scale=1]{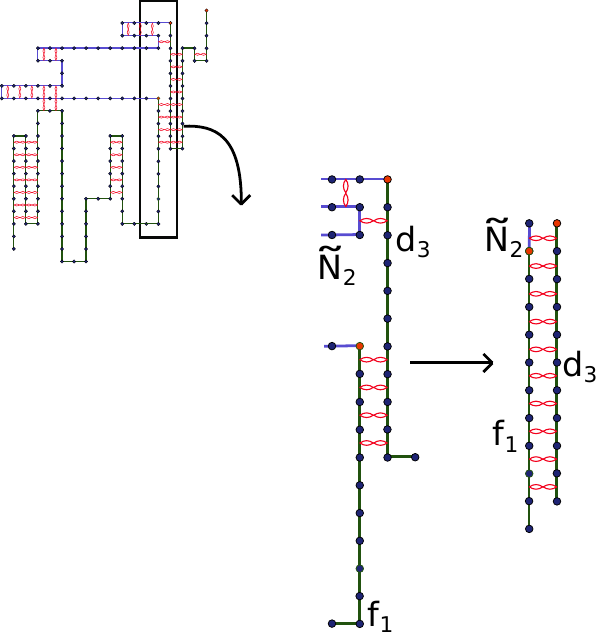}
\caption{On the left, a NE-PSAW path made of three blocks. 
In the picture we zoom in on the interactions between the third block and the 
rest of path. 
We recall that the third block can only interact with its two preceding blocks, i.e., the first and the second one.
We call $f_1$ the last vertical stretch of the first block and
$d_3$ the first vertical stretch of the third block.
% $\tilde N_2$ the fraction of $N_2$ which may truly interact with $d_3$.
 The interactions between the first and the third blocks involve   $f_1$ and $d_3$ while the interactions 
 between the second and the third blocks involve $d_3$ and  $\tilde N_2$ (the number of inter-stretches in the second block that 
 may truly interact with $d_3$, on the picture $\tilde N_2=1$).
Such interactions are bounded above by $(\tilde N_2 + f_{1}) \tilde{\wedge} d_{3}$. 
}
\label{fig:IPDRW3}
\end{figure}

\subsection{Conclusion} 
 
From \eqref{include} it is straightforward that  $F^{\rm\tx{IPSAW}}(\beta)\geq f(\beta)$ for every $\beta\geq 0$ (recall \eqref{deff}). Thus,  Theorem \ref{Thm1}  implies  that the critical point of IPSAW is not smaller than that of IPDSAW, i.e., 
\begin{equation}
\label{eq:IPIPDcriticalpoint}
\beta_c \, \leq\, \beta_c^{\tx{IPSAW}}.
\end{equation}
It would be interesting to understand whether \eqref{eq:IPIPDcriticalpoint} is an equality or not. Even more challenging would be the computation of $\beta_c^{\tx{IPSAW}}$.

Let us stress also that even if the existence of a collapsed transition for IPSAW is proven, we do not have any results concerning its scaling limit in each regime (extended, critical  and collapsed) as we did for IPDSAW in Section \ref{sec:scalingIPDSAW}. 
In this spirit, at $\beta=0$, the scaling limit of the prudent walk itself has been derived in  \cite{PTS17}. We conjecture that 
in the extended phase the scaling limit should have a similar  structure, that is, a straight line. More interesting is the inside of the collapsed phase, in which the limit shape is less clear. In analogy with the results obtained in Section \ref{sec:scalingIPDSAW} we only expect it to be deterministic.

 \medskip

We  conclude with a few words about the $2$-dimensional Interacting Self-Avoiding Walk (ISAW) defined exactly like the IPSAW in (\ref{eq:ham1}) but with a larger set of allowed configurations, that is 
\begin{align}\label{allconf}
\Omega_L^{\tx{SAW}}:=\big\{w:=(w_i)_{i=0}^L\in (\mathbb{Z}^2)^{L+1}\colon\,& w_{0}=0,\  w_{i+1}-w_i\in \{\leftarrow,\rightarrow,\downarrow,\uparrow\},\  0\leq i \leq L-1,\\
\nonumber &w \  \text{satisfies the self-avoiding condition}\big\}.
\end{align}
We denote by $Z_{L,\beta}^{\tx{ISAW}}$ the partition function of ISAW  and we define its  free energy as
\begin{equation}\label{defISAW}
F^{\tx{ISAW}}(\beta):=\liminf_{L\to \infty} \frac{1}{L} \log Z_{L,\beta}^{\tx{ISAW}},
\end{equation}
%where the $\liminf$ in \eqref{defISAW} is chosen to overstep the fact that 
%the convergence of the free energy remains an open issue.
Theorem \ref{Thm33} below shows that the  conjectured collapse transition displayed by ISAW at some 
$\beta_c^{\tx{ISAW}}$ does 
not correspond to a self-touching saturation as is the case for 
IPDSAW and IPSAW. The reason is that, even very dense ISAW trajectories can integrate small holes 
which are not compatible with the prudent condition. Introducing a small density of holes ($\gep L$) in a dense
ISAW configuration of length $L$   yields a loss
of self-touchings of order $\gep L$, however this is overcome by the entropy gain associated with the choice of 
the locations of those holes (of order $-\gep \log(\gep) L$).  
\begin{theorem}[\cite{PT17}, Theorem 2.3]
\label{Thm33}
\begin{align}
\mathrm  F^{\tx{ISAW}}(\beta)&>\beta, \quad  \text{for every} \quad  \beta\in [0,\infty).
\end{align} 
\end{theorem}

\bibliographystyle{imsart-nameyear}
\bibliography{cnp}

%%%%%%%%%%%%%%%%%%%%%%%%%%%%%%%%%%%%%%%%%%%%%%%%%%%%%%%%%%%%%%%%%%%%%%%%%%

\end{document}